\documentclass[trsc,nonblindrev]{informs3}
\renewcommand{\theARTICLETOP}{}
\renewcommand{\theRRHSecondLine}{}
\renewcommand{\theLRHSecondLine}{}
\usepackage{amsmath}

\OneAndAHalfSpacedXI %

\usepackage{endnotes}
\usepackage{bm}

\usepackage{color-edits}
\addauthor[Ang]{ax}{magenta}
\addauthor[Chiwei]{cy}{cyan}

\usepackage{algorithm}
\usepackage{algpseudocode}
\usepackage{tikz}
\usetikzlibrary{tikzmark}
\usepackage{subcaption}
\usepackage{diagbox}
\usetikzlibrary{arrows}
\usepackage{booktabs}
\usepackage{comment}
\usepackage{enumitem} 
\usepackage{pgfplots}
\usepgfplotslibrary{groupplots}
\pgfplotsset{compat=1.18}
\usepackage{xcolor}

\definecolor{mplblue}{HTML}{1F77B4}
\definecolor{mplorange}{HTML}{FF7F0E}
\definecolor{mplgreen}{HTML}{2CA02C}

\usepackage{natbib}
 \bibpunct[, ]{(}{)}{,}{a}{}{,}%
 \def\bibfont{\small}%

\EquationsNumberedThrough    %

\TheoremsNumberedThrough     %
\ECRepeatTheorems  %

\begin{document}

\RUNAUTHOR{Xu and Yan}

\RUNTITLE{Dispatching and Pricing in Spatial Queues}

\TITLE{Dispatching and Pricing in Two-Sided Spatial Queues}

\ARTICLEAUTHORS{%
\AUTHOR{Ang Xu}
\AFF{Department of Industrial Engineering and Operations Research,
University of California, Berkeley, \EMAIL{angxu@berkeley.edu}}

\AUTHOR{Chiwei Yan}
\AFF{Department of Industrial Engineering and Operations Research,
University of California, Berkeley, \EMAIL{chiwei@berkeley.edu}}
} %

\ABSTRACT{%
We study a dispatching and pricing problem in two-sided spatial queues with fixed supply, motivated by ride-hailing and robotaxi platforms. Idle drivers queue on one side, waiting to pick up riders, while riders queue on the other, waiting to be matched with available drivers. The platform seeks to maximize net profit, penalized by rider waiting penalties, by jointly optimizing state-dependent dispatching and pricing decisions. We formulate this problem as a Markov decision process with state-dependent service times that capture key features of spatial matching. We show that, under mild assumptions, the optimal dispatching policy admits a closed-form expression with a \emph{zigzag} structure. This policy significantly improves the tractability of pricing optimization due to the resulting closed-form stationary distribution and a substantially reduced state space. Building on this insight, we propose an efficient and scalable dynamic programming heuristic to approximate the optimal zigzag policy in more general settings. 
Extensive numerical experiments with both the analytical model and ride-hailing simulations demonstrate that our algorithm is both near-optimal and highly scalable.

}%

\HISTORY{First version, 07/2025. Revision, 03/2026, 05/2026.}

\maketitle

\newcommand{\marker}[1]{}

\section{Introduction}\label{sec:Intro}

Ride-hailing services, together with the rapidly growing robotaxi industry, have revolutionized urban transportation. The dispatching process, which assigns available drivers to incoming ride requests, as well as the pricing algorithms for each incoming rider are both critical to the operational efficiency of these platforms~\citep{yan2020dynamic}. Ineffective dispatching and pricing can result in the so-called ``wild goose chases''  phenomenon~\citep{castillo2025matching}, where drivers experience long pick-up times and low utilization (the proportion of on-trip time relative to total supply minutes) due to supply shortages caused by overly aggressive dispatching (e.g., assigning drivers to riders despite long pick-up times) or insufficient price adjustments to moderate demand during peak periods.

In this paper, we study the joint dispatching and pricing problem in a two-sided spatial queueing model, where both idle drivers and riders may queue in the system while awaiting dispatch. In contrast to ``one-sided'' queueing models---where riders and idle drivers are matched immediately and either the rider or driver side queues (e.g., \citealt{feng2020we,besbes2022spatial,benjaafar2024human,yan2026matching})---we allow both sides to queue simultaneously, and enable the platform to intentionally idle some drivers even when there are riders waiting, a strategy used often in practice when the pickup distance is long (see, e.g.,  \citealt{xu2020supply, wang2022ondemand,castillo2025matching}). We briefly summarize our model.

\smallskip

\noindent\textbf{Model.} We develop a continuous-time Markov decision process to approximate the dispatching and pricing decisions of a ride-hailing platform with a \emph{fixed} fleet size (e.g., robotaxi platforms). The system state tracks the number of drivers in service, equivalently the number of idle drivers, and the number of waiting riders. The platform dynamically makes dispatching decisions and adjusts ride prices based on the evolving system state, aiming to maximize long-run average net profit minus rider waiting penalties. A key feature of our model is that the service rate is \emph{state-dependent}, capturing the dynamics of spatial matching. Intuitively, as the number of idle drivers or waiting riders increases, the expected pick-up time decreases. We incorporate this dependence by modeling the service (trip completion) rate as a function of the idle driver and waiting rider counts in the system. 

\smallskip

\noindent\textbf{Our result.} We show that, under mild assumptions on the service rate and cost coefficients, the optimal dispatching policy admits a closed-form solution (Theorem~\ref{theorem:opt_policy_struct}). Inspired by the structure of this solution, we introduce a class of dispatching policies termed \emph{zigzag} policies. As the system state evolves, the platform’s decisions alternate between ``hold'' and ``dispatch'' in a pattern that, when visualized over the state space of drivers and riders, traces out a zigzag curve. Moreover, the zigzag policy exhibits a monotone, threshold-type structure, making it interpretable. This policy also significantly improves the tractability of pricing optimization due to the resulting closed-form stationary distribution and a substantially reduced (recurrent) state space. Building on this insight, we propose an efficient and scalable dynamic programming heuristic to approximate the optimal zigzag policy
in more general settings. We evaluate our algorithm through both analytical model-based experiments and a ride-hailing simulation that captures realistic dispatching processes. In all cases, benchmarking against competitive baselines demonstrates that the zigzag policy achieves near-optimal performance and exhibits excellent scalability.

\smallskip

\noindent\textbf{Organization.} The remainder of the paper is organized as follows. Section~\ref{sec:lit_review} reviews related literature. Section~\ref{sec:model_and_prelim} introduces the model, and Section~\ref{sec:methodology} presents our analytical results and solution approaches. Section~\ref{sec:numerical_experiments} reports extensive numerical experiments using both the analytical model and ride-hailing simulations. We conclude in Section~\ref{sec:conclusion}.

\smallskip

\section{Literature Review}
\label{sec:lit_review}

Our research is closely related to the topics of dynamic spatial matching and pricing for reusable resources, many of which are motivated by ride-hailing and shared mobility operations. We summarize the most relevant works on these two topics in the following paragraphs. %

\textbf{Dynamic spatial matching.} 
Our work belongs to the stream of literature that uses queueing-theoretic models to approximate a dynamic spatial matching system motivated from ride-hailing platforms. A prevailing modeling choice is to use state-dependent service time to capture spatial friction (e.g., \citealt{
feng2020we,besbes2022spatial,benjaafar2024human%
})---the service rate of each trip increases in the number of idle drivers or the number of riders waiting in the queue, since more idle drivers or waiting riders imply a shorter pickup time, and thus, service time. Many prior studies focus on performance evaluation or capacity planning under \emph{non-idling} dispatching policies, where idle drivers are immediately matched with waiting riders. In such settings, the queues are effectively one-sided---either on the rider or driver side, depending on availability. We consider a more general setting in which the platform is allowed to intentionally idle drivers even when riders are waiting, resulting in two-sided queues. This has connection to prior work through the use of \emph{dispatching radius caps} (e.g., \citealt{feng2020we,wang2022ondemand,castillo2025matching}) or \emph{temporal pooling} \citep{chen2025depooling}, via heuristic or fluid analysis, or under ``one-sided'' settings (see, e.g., \citealt{xu2020supply} and \citealt{chen2025depooling}). In contrast, our work aims to characterize the structure of the optimal (potentially idling) dispatching policy within a Markov decision process framework that captures two-sided queues with state-dependent service times---a formulation which is a natural extension but, to our knowledge, has not yet been addressed in the existing literature.

An important assumption of the aforementioned spatial queueing models (including our work) is that service times are independently distributed. In real dynamic spatial matching systems, however, there is often significant correlation across dispatches, as the locations of idle drivers and waiting riders tend to change little over short time intervals. Focusing on spatial correlation and non-idling dispatching policy, \cite{yash2022dynamic} considers a different framework that explicitly models the spatial locations of riders and drivers, and proposes a rate-optimal dispatching policy relative to benchmark approaches or performance bounds. On the other hand, the temporal aspect is simplified in \cite{yash2022dynamic}, where all trips are assumed to be completed within a single time period. In a related spatial setting on a line, \cite{balkanski2025powergreedyonlineminimum} show that a greedy policy achieves a constant competitive ratio relative to the offline optimum when each driver exits the system after serving a customer. Although our model simplifies spatial complexity through ``reduced-form'' Markovian queueing dynamics, we evaluate the proposed policy in a simulation environment with fully detailed spatial features and demonstrate its effectiveness and robustness.

\smallskip

\textbf{Pricing in (spatial) service systems.} The pricing aspect of our work also relates to the literature on pricing in (spatial) service systems. Conceptually speaking, our work is closely related to \cite{castillo2025matching}, which demonstrates the need for surge pricing to avoid ``wild goose chases'' and ensure low pick-up times, thereby restoring the demand-supply balance. Methodologically, our work is most relevant to \cite{banerjee2015pricing}, \cite{besbes2022static}, \cite{elmachtoub2025power}, and \cite{bergquist2024staticpricing}, which study pricing for reusable resources modeled using Markovian queueing dynamics. Many of these studies focus on analyzing the performance of static pricing policies. A distinguishing feature of our work is the incorporation of state-dependent service rates, in contrast to prior studies that assume a fixed service rate, rendering many of their insights not directly applicable. State-dependent service rates also make the service policy nontrivial, making the pricing policy tightly coupled with the complexity of the service policy. One implication of our proposed zigzag dispatching policy is that the induced Markov chain effectively reduces to a birth-death process with a state space that is an order of magnitude smaller than the original. This structure also admits a closed-form stationary distribution, making the pricing optimization problem more tractable.

\section{Model and Preliminaries}
\label{sec:model_and_prelim}

We consider a ride-hailing or robotaxi platform operating a fixed fleet to serve riders across a service region. Rider arrivals into the service region follow a Poisson process with rate $\Lambda$. Upon arrival, each rider is quoted a price that may depend on the current system state. If a rider accepts the price, she joins a queue to await matching with a driver. The platform can decline to serve a rider by setting the price at or above a certain threshold. Importantly, the platform is not required to dispatch a driver immediately when a rider agrees to the price; it may delay dispatch even when idle drivers are available. We assume riders are patient---once they join the queue, they will wait until matched with an idle driver. However, the platform incurs a penalty for rider waiting times in its objective, and thus seeks pricing and dispatching policies that maximize net profit minus rider waiting penalties, as detailed later.

\subsection{State, Action, and Policy}\label{subsec:state,action, and policy} %

Let $L$ denote the total number of drivers in the system, and let $\mathcal{L} := \{0, 1, \ldots, L\}$ be the index set of possible driver counts. %
The state space is defined as $\mathcal{S} := \{(l, m) : l \in \mathcal{L},\, m \in \mathbb Z_{\ge0}\}$, where each state $(l, m)$ represents $l$ drivers currently \emph{in service} and $m$ riders \emph{waiting in the queue}. State transitions occur due to three types of events: (1) \emph{Trip completion:} when a driver completes a trip, he immediately becomes idle and the state transitions from $(l, m)$ to $(l-1, m)$; (2) \emph{Rider arrival:} when a new rider arrives and accepts the price to join the queue, the state transitions from $(l, m)$ to $(l, m+1)$; (3) \emph{Dispatching decision:} when the platform dispatches $k$ drivers to pick up $k$ waiting riders (with $k \leq m$), the state transitions from $(l, m)$ to $(l + k, m - k)$. The platform may dispatch \emph{multiple} drivers \emph{simultaneously} at any decision epoch. For example, if the system is in state $(l, m)$, a trip is then completed, and the platform decides to immediately dispatch two drivers to pick up two riders, the state transitions to $(l - 1 + 2, m - 2) = (l + 1, m - 2)$.

The platform's actions consist of two components: \emph{pricing} and \emph{dispatching}. At each state transition, the platform updates both the price for incoming riders and its dispatching decisions. We assume the trip price is given by $p = p_0 + p_1 d$, with $p_0 \geq 0$ denoting a fixed base fare, and $p_1 \in [0, p_{\mathrm{max}}]$ representing the per-kilometer charge that can be dynamically adjusted by the platform (e.g., via surge pricing). Riders' willingness to pay per kilometer is drawn from a cumulative distribution function (CDF) $F(\cdot)$. A rider will accept the price and join the queue if and only if their willingness to pay per kilometer exceeds or equals $p_1$. Thus, the \emph{effective request rate} (or request conversion rate) is $\lambda = \Lambda \bar{F}(p_1)$,\footnote{The dependence of demand on $p_0$ is captured by the exogenous arrival rate $\Lambda$, as the model does not permit $p_0$ to be altered.} where $\bar{F}(p_1) = 1 - F(p_1)$ denotes the tail CDF of riders' willingness to pay. We assume a one-to-one correspondence between $\lambda$ and $p_1$, with $\bar{F}(0) = 1$, $\bar{F}(p_{\mathrm{max}}) = 0$, and $\bar{F}(p_1)$ strictly decreasing in $p_1$. Therefore, setting $p_1$ is equivalent to setting $\lambda$, and we denote by $p_1(\lambda)$ the per-kilometer price corresponding to an effective request rate $\lambda$. For notational convenience, we define the \emph{pricing policy} as a mapping $\lambda: \mathcal{S} \to [0, \Lambda]$, where $\lambda(l, m)$ specifies the effective request rate at state $(l, m)$. In the following, we describe pricing decisions in terms of $\lambda$ rather than $p$.

A \textit{dispatching policy} $\phi: \mathcal{S} \mapsto \{0,1\}$ maps the state space to a binary number, which determines whether the platform dispatches available drivers to riders at state $(l,m)$.
In specific, $\phi(l,m)=1$ represents the dispatching action, and $\phi(l,m)=0$ the holding action (no dispatching) under state $(l,m)$.
The system transitions from state $(l, m)$ to state $(l+1, m-1)$ if $\phi(l,m)=1$. %
At the new state $(l+1, m-1)$, the platform considers again whether to dispatch and repeats this process until policy $\phi$ at the current state indicates a holding action. For example, suppose the dispatching policy $\phi$ is defined as $\phi(1,3) = \phi(2,2) = 1$ and $\phi(3,1) = 0$. If a new rider arrives when the system is in state $(l, m) = (1,2)$, the state first changes to $(1,3)$. Since $\phi(1,3) = 1$, a driver is dispatched, causing the state to transition to $(2,2)$. Similarly, because $\phi(2,2) = 1$, another driver is dispatched, and the state transitions to $(3,1)$. Note that the platform is able to dispatch only when there is at least one idle driver and one rider in the queue. Thus, for any policy, we have $\phi(l,0) = 0$ for all $l \in \mathcal{L}$ and $\phi(L,m) = 0$ for all $m \in \mathbb Z_{\ge0}$. Note that the system will never remain in a state $(l,m)$ with $\phi(l,m) = 1$, as all dispatching actions are executed instantaneously without delay.

\smallskip

\subsection{Service (Trip Completion) Rate} 
For each trip, the service time consists of two parts: ``pick-up time'' and ``on-trip time''. Pick-up time is the duration it takes for the driver to travel from his current location to the rider's arrival location. On-trip time is the duration for the driver to travel from the rider's arrival location to her desired destination. Instead of modeling the real-time locations of riders and drivers, following related literature such as \cite{feng2020we},  \cite{besbes2022spatial} and \cite{benjaafar2024human}, 
we adopt a ``reduced-form'' setting which approximates the service time of each driver in service at state $(l,m)$ as an exponential\footnote{In the real world, service times may be correlated with past history of arrivals and dispatches. For tractability, we ignore this dependence and assume that service times are independent and memoryless.} random variable with rate $\mu_{l,m}$ which increases with $m$ and decreases with $l$.\footnote{In fact, the monotonicity of $\mu_{l,m}$ is \emph{not} required in our analysis, although it is well expected in real-world spatial systems.} %
Intuitively, as the number of idle drivers or waiting riders increases, the expected pick-up time decreases. For example, if the platform matches the closest idle driver and rider pair when making a dispatch decision---a common practice in the industry \citep{yan2020dynamic}---a larger pool of waiting drivers or riders enables better matches. %
We also assume that there exists a $\bar{\mu}>0$ such that $\mu_{l,m}\leq \bar{\mu}$ for all $(l,m)\in\mathcal{S}$. 
The boundedness assumption for $\mu_{l,m}$ is motivated by the fact that the service time necessarily includes a non-negligible on-trip component. The trip completion rate at state $(l,m)$ is thus $l\mu_{l,m}$. We assume that the on-trip time for each ride has a mean $t_0$, and the driving speed of all vehicles is $1$. In general, for all $(l,m)\in\mathcal{S}$, we have $t_0 \le 1/\mu_{l,m}$, since the on-trip time is one component of the total service time.

\subsection{Objective}\label{sec:objective}
Our objective balances two components: (1) improving the platform’s expected long-run average profit (i.e., revenue minus operating costs), and (2) reducing penalties associated with rider queueing and pickup delays.

\smallskip

\noindent\textbf{Platform revenue.}  
Let $\mathcal{R}(\lambda,\phi,T)$ be the cumulative revenue up to time $T$ under policy $(\lambda,\phi)$. We denote the corresponding long-term average revenue rate by $\mathcal{R}(\lambda,\phi):=\lim_{T\rightarrow\infty} \mathbb{E}\left[\mathcal{R}(\lambda,\phi,T)\right]/T$. Note that some policies may lead to unstable queues, where the number of riders in the queue grows indefinitely. However, without loss of optimality, we can restrict our analysis to policies that ensure queue stability, as unstable policy leads to infinitely large rider waiting penalties, as we will introduce next. When the queue is stable, the stationary distribution exists as the underlying unichain is ergodic. By the ergodic theorem, the long-run time-average revenue converges almost surely to the expected revenue under the stationary distribution. As a consequence, $\mathcal{R}(\lambda, \phi)$ is well-defined and finite. Let $\pi_{\lambda,\phi}(l,m)$ be the stationary probability of having $l$ drivers in service and $m$ riders in the queue under policy $(\lambda,\phi)$. Then we can write $\mathcal{R}(\lambda,\phi)$ as
\begin{equation}\label{eq:rev}
\mathcal{R}(\lambda,\phi)=\sum_{(l,m)\in\mathcal{S}}\pi_{\lambda,\phi}(l,m)\lambda(l,m)\left(p_0+p_1(\lambda(l,m))\cdot t_0\right).
\end{equation}

\noindent\textbf{Operating cost and rider waiting penalties.} %
We consider the following quantities under policy $(\lambda, \phi)$: (1) the expected number of drivers in service (including those en route to pick-ups and on trips), denoted by $L^{d}_s(\lambda, \phi)$; (2) the expected number of idle (available) drivers, denoted by $L^{d}_o(\lambda, \phi)$, where $L^{d}_o(\lambda, \phi) = L - L^{d}_s(\lambda, \phi)$; (3) the expected number of riders waiting to be picked up, denoted by $L^r_p(\lambda, \phi)$; and (4) the expected number of riders in the queue awaiting a match, denoted by $L^r_q(\lambda, \phi)$. Our objective %
can be expressed as
\begin{equation}\label{eq:objective}
\sup_{\lambda,\phi}\quad\underbrace{\mathcal{R}(\lambda,\phi)}_{\text{revenue rate}}-\underbrace{w_s^d L_s^d(\lambda,\phi) - w_o^d L_o^d(\lambda,\phi)}_{\text{driving/operating cost rate}}-\underbrace{w_p^r L_p^r(\lambda,\phi) - w_q^r L_q^r(\lambda,\phi)}_{\text{penalty rate from rider waiting}},  %
\end{equation}
where $w_s^d, w_o^d, w_p^r$ and $w_q^r$ are positive constants measuring cost/penalties per unit time and driver/rider. %
In particular, $w_s^d L_s^d(\lambda,\phi)$ and $w_o^d L_o^d(\lambda,\phi)$ denote the operating costs for drivers while in service and while idle, respectively, with $w_s^d \geq w_o^d$ typically holding since idle drivers generally consume less fuel or electricity. Similarly, $L_p^r(\lambda, \phi)$ and $L_q^r(\lambda, \phi)$ capture the penalties associated with riders waiting to be picked up and those waiting in the queue, respectively. Usually, $w_q^r \geq w_p^r$, as riders are generally more sensitive to waiting in the queue than after a driver has been dispatched to them.
Calculating the objective in~\eqref{eq:objective} directly can be challenging, as the expected number of riders waiting to be picked up, $L_p^r(\lambda, \phi)$, cannot be explicitly determined from the system state. Fortunately, an equivalent and much simpler form of the objective function can be derived. Let $\mathcal{R}(\lambda, \phi; p_0)$ denote the value of $\mathcal{R}(\lambda, \phi)$ for a given $p_0$.

\begin{lemma}\label{lemma:equivalent_penalty}
    We have $\argmax_{\lambda, \phi}\mathcal{R}(\lambda,\phi;p_0)-w_s^d L_s^d(\lambda,\phi) - w_o^d L_o^d(\lambda,\phi) - w_p^r L_p^r(\lambda,\phi) - w_q^r L_q^r(\lambda,\phi) = \argmax_{\lambda, \phi}\mathcal{R}(\lambda,\phi; p_0+w_p^r t_0)-(w_s^d+w_p^r-w_o^d)L_s^d(\lambda,\phi)-w_q^r L_q^r(\lambda,\phi)$.
\end{lemma}

Lemma~\ref{lemma:equivalent_penalty} shows that, by appropriately resetting the value of $p_0$ and the cost/penalty coefficients, both $L_o^d(\lambda, \phi)$ and $L_p^r(\lambda, \phi)$ can be omitted from the objective, while the optimal policy remains unchanged. %
To see the intuition, since drivers are either in service or idle, in \eqref{eq:objective}, the driving/operating cost rate can be rewritten as the following two parts: (i) a penalty term $w_o^dL$ that accrues for all drivers and (ii) a penalty term $(w_s^d-w_o^d)L_s^d(\lambda,\phi)$ that accrues only for drivers who are in service. Without loss of optimality, we can ignore the penalty term in part (i), as it is independent of $\lambda$ and $\phi$. Thus, $L_o^d(\lambda, \phi)$ can be omitted from our objective. Likewise, when a driver is in service, she is either picking up a rider or on a trip to the destination. Therefore, we can rewrite $L_p^r(\lambda,\phi)$ as $L_s^d(\lambda,\phi)$ minus the average number of on-trip drivers. One can show that the on-trip component is equivalent to increasing $p_0$ by $w_p^rt_0$. This allows us to omit $L_p^r(\lambda, \phi)$ as well. %
In the resulting optimization problem, the expected number of drivers in service $L_s^d(\lambda,\phi)$ and the expected number of riders in the queue $L_q^r(\lambda,\phi)$ can be directly calculated from the states $(l,m)$ and the corresponding stationary probability. For conciseness, in the remainder of the paper, we set driver-side penalty $c_d:=w_s^d + w_p^r - w_o^d$ and rider-side penalty $c_r:=w_q^r$. %
Our new (equivalent) objective, denoted by $\tilde{\mathcal{R}}(\lambda, \phi)$, becomes %
\begin{align}
\label{eq:objective_alt}
\sup_{\lambda,\phi}\quad\tilde{\mathcal{R}}(\lambda, \phi) 
    &:= \sum_{(l,m)\in\mathcal{S}} \pi_{\lambda, \phi}(l,m)\left( \lambda(l,m) \left(p_0 + w_p^rt_0 + p_1(\lambda(l,m))\cdot t_0\right)- c_d l - c_r m \right).
\end{align}

Our goal is to develop efficient solution approaches to find an optimal policy $(\lambda^\ast,\phi^\ast)$ that maximizes this objective, as we will discuss next. %

\section{Solution Approaches}
\label{sec:methodology}
The system described in Section~\ref{sec:model_and_prelim} is a \emph{continuous-time Markov decision process} (CTMDP) with a countably infinite state space. However, when $c_r > 0$, it is without loss of optimality to impose an upper bound $M$ on the queue length and reject all incoming riders once the queue reaches this limit. Otherwise, allowing the queue to grow indefinitely would result in unbounded penalties in the objective. 
To find an optimal policy in a CTMDP, as outlined in \cite{Puterman1994MDP}, a common approach is to transform it into a discrete-time MDP by using the \emph{uniformization} method. %
However, to apply this method, the maximum transition rate among all states has to be bounded. Unfortunately, under the action space defined in Section \ref{sec:model_and_prelim}, the transition rate can become \emph{infinite} as multiple state transitions can occur instantaneously---e.g., if an arrival or trip completion is followed by (possibly multiple) dispatching actions, then the transition rates at these intermediate dispatching states are effectively infinite, since these dispatching actions occur with no time delay.
To bound the transition rates at each state, we have to expand the action space. Upon each event, either a new rider request arrival or a trip completion, the platform makes a dispatch decision. Specifically, at each state, it must now determine \emph{how many} drivers to dispatch in response to a rider arrival, and separately, \emph{how many} to dispatch upon a trip completion.
The formal definition of this new action space and the corresponding Bellman optimality equation are detailed in Appendix \ref{subsec:alt_def_action_space}.%
\footnote{Note that this new action space contains the action space defined in Section \ref{subsec:state,action, and policy}, nevertheless it can be shown that there always exists an optimal policy that is contained in both spaces. Consequently, they are equivalent from an optimization standpoint.} %
Although this alternative action space guarantees bounded transition rates, it introduces significant redundancy: instead of a binary dispatch decision for each state, the action can now take any integer value up to the minimum of available riders and drivers, depending on future events. This added complexity makes it intractable. While one could restrict the action space to dispatching at most one driver at a time, we will show numerically in Section~\ref{subsec:synthetic_exp} that the problem remains computationally intensive for real-world instances. This motivates the need for alternative, efficient solution approaches, which we develop by leveraging the structural properties of our problem.

\subsection{Properties of the Optimal Dispatching Policies}\label{sec:optimal_dispatching_policies}

We begin by analyzing the properties of optimal dispatching policies. Notably, given a dispatching policy, the pricing optimization problem can become simpler: it suffices to consider only the class of recurrent states under the specified dispatching policy~$\phi$, which can be significantly smaller than the full state space (in fact, by an order of magnitude smaller under our later proposed dispatching policy). We first introduce a mild assumption on $\mu$.

\begin{assumption}\label{as:2}
We assume that $\mu$ satisfies the following two conditions.
\begin{itemize}
\item[(1)] $l\mu_{l, m+1}-l\mu_{l, m}$ is non-increasing in $m$ and non-decreasing in $l$;
\item[(2)] $(l+1)\mu_{l+1, m}-l\mu_{l, m}$ is non-increasing in $l$ and non-decreasing in $m$.\footnote{In fact, $l\mu_{l, m+1}-l\mu_{l, m}$ being non-decreasing in $l$ in condition (1) is equivalent to $(l+1)\mu_{l+1, m}-l\mu_{l, m}$ being non-decreasing in $m$ in condition (2). We keep both in the assumption for clarity.}
\end{itemize}
\end{assumption}

Roughly speaking, Assumption~\ref{as:2} captures the idea of diminishing returns: adding an additional idle driver or waiting rider yields progressively smaller increments in service rate. In fact, Assumption~\ref{as:2} can be more easily understood by considering analogous monotonicity conditions on the differences of $\mu_{l,m}$, rather than $l\mu_{l,m}$, which is a \emph{stronger} requirement that implies Assumption~\ref{as:2}. 
For condition (1) on $m$, as rider queue length increases, the marginal impact of an additional waiting rider on pickup distance and service rate %
decreases. This is equivalent to the discrete concavity of $l\mu_{l,m}$ in $m$ for fixed $l$.
For condition (1) on $l$, as the number of drivers in service $l$ increases (and idle drivers $L - l$ decrease), the marginal impact of an additional rider on service rate, $\mu_{l,m+1} - \mu_{l,m}$, increases, so as $l\mu_{l,m+1} - l\mu_{l,m}$. %
For condition (2), we focus on $l$ since monotonicity in $m$ mirrors that in $l$ under condition (1). As $l$ increases, the marginal impact of an additional idle driver, $\mu_{l,m} - \mu_{l+1,m}$, increases, implying that $(l+1)\mu_{l+1,m} - l\mu_{l,m}$ decreases in $l$. This corresponds to the discrete concavity of $l\mu_{l,m}$ in $l$ for fixed $m$.\footnote{Assumption~\ref{as:2} may be compared to the joint discrete concavity of $l\mu_{l,m}$ in both $l$ and $m$. However, neither condition implies the other: joint concavity requires the Hessian matrix (formed by second forward differences) of $l\mu_{l,m}$ to be negative semi-definite, whereas Assumption~\ref{as:2} is equivalent to requiring non-positive diagonal entries and non-negative off-diagonal entries in this Hessian.}

We now proceed to discuss how we can derive an optimal dispatching policy under Assumption \ref{as:2}. We first classify a state into either type 1 or type 2, defined as follows. 

\smallskip

\begin{definition}[State Classification]
     We classify a state $(l, m) \in \mathcal{S}$ as \emph{type 1} if its trip completion rate, $l\mu_{l,m}$, exceeds that at $(l+1, m-1)$, i.e., $l\mu_{l,m} > (l+1)\mu_{l+1, m-1}$. If $l=L$ or $m=0$, $(l,m)$ is classified as a type 1 state as well. Otherwise, if $l\mu_{l,m} \leq (l+1)\mu_{l+1, m-1}$, we classify $(l, m)$ as \emph{type 2}.
\end{definition}

\smallskip

In short, a type 1 state has a higher trip completion rate than its lower-left neighbor, the state reached by dispatching exactly one driver. Conversely, a type 2 state has a lower trip completion rate than its lower-left neighbor. The following lemma provides an important observation about the distribution of these two types of states.

\begin{lemma}\label{lemma:serv_rate_ineq}
 Suppose that $\mu$ satisfies Assumption \ref{as:2}. Then the following statements hold.
 \begin{itemize}
     \item[(1)] If $(l',m')$ is a type 1 state, then for any $l\geq l'$ and $m\leq m'$, $(l,m)$ is also a type 1 state.
     \item[(2)] If $(l',m')$ is a type 2 state, then for any $l\leq l'$ and $m\geq m'$, $(l,m)$ is also a type 2 state.  
 \end{itemize}
\end{lemma}

\begin{table}[htbp]
     \small
    \renewcommand{\arraystretch}{1.5} %

    \begin{subtable}{.5\linewidth}
      \centering

        \begin{tabular}{c|ccccc}
         \backslashbox[12mm]{$l$}{$m$} & %
         & $m'-1$ & $m'$ & \\
         \hline %
         \\
         $l'$ & $\cdots$ & $l'\mu_{l',m'-1}$ & $l'\mu_{l',m'}$ &  \\
         $l'+1$ & $\cdots$ &  $(l'+1)\mu_{l'+1,m'-1}$ & $(l'+1)\mu_{l'+1,m'}$  &  \\
         & &  $\vdots$ &  $\vdots$\\
        \end{tabular}

        \begin{tikzpicture}[overlay, remember picture]
        \node[text=black, rotate=35] at (0.9,1.3){$\pmb{<}$};
        \node[rotate=35] at (-2,1.3){$<$};
        \node[rotate=35] at (-2,0.5){$<$};
        \node[rotate=35] at (0.9,0.5){$<$};
        \end{tikzpicture}
        \caption{%
       \small States to the left of or below a type 1 are type 1.}
        \label{subtab:def_type1_state}
    \end{subtable}%
       \begin{subtable}{.5\linewidth}
      \centering
        \begin{tabular}{c|ccccc}
         \backslashbox[10mm]{$l$}{$m$} & \qquad \qquad & $m'-1$ & $m'$ & %
         \\
         \hline %
         & & $\vdots$ & $\vdots$ \\
         $l'$ &  & $l'\mu_{l',m'-1}$ & $l'\mu_{l',m'}$ & $\cdots$ \\
         $l'+1$ & &  $(l'+1)\mu_{l'+1,m'-1}$ & $(l'+1)\mu_{l'+1,m'}$ & $\cdots$ \\
         \\
        \end{tabular}

        \begin{tikzpicture}[overlay, remember picture]
        \node[text=black, rotate=30] at (0.9,1.3){$\pmb{\geqslant}$};
        \node[rotate=30] at (3,1.3){$\geqslant$};
        \node[rotate=30] at (0.9,1.8){$\geqslant$};
        \node[rotate=30] at (3,1.8){$\geqslant$};
        \end{tikzpicture}
        \caption{%
        \small States to the right of or above a type 2 are type 2.}       \label{subtab:def_type2_state}
    \end{subtable}
\caption{Illustration of Lemma \ref{lemma:serv_rate_ineq}. \normalfont{If the inequality shown in bold is satisfied, then the nonbold inequalities are satisfied as well.}}
\end{table}

Lemma~\ref{lemma:serv_rate_ineq} implies that under Assumption~\ref{as:2}, each anti-diagonal of the state map contains a state whose trip completion rate is greater than or equal to that of all other states along the same anti-diagonal. Intuitively, staying at states with high trip completion rates might benefit system performance because trips are completed more quickly, leading to faster driver turnover. Our main result below formalizes this.

\begin{theorem}[Optimal Dispatching Policy]\label{theorem:opt_policy_struct}
     Suppose that %
     $\mu$ satisfies Assumption~\ref{as:2}. Then, the following statements hold. 
     \begin{enumerate}
         \item[(1)] If $c_d=c_r$, there exists an optimal policy $(\lambda^\ast, \phi^\ast)$ such that $\phi^\ast (l,m)=0$ for all type 1 states and $\phi^\ast (l,m)=1$ for all type 2 states. %
         \item[(2)] If $c_d>c_r$, %
         there exists an optimal policy $(\lambda^\ast, \phi^\ast)$ such that $\phi^\ast (l,m)=0$ for all type 1 states.
         \item[(3)] If $c_d<c_r$, %
         there exists an optimal policy $(\lambda^\ast, {\phi}^\ast)$ such that $\phi^\ast (l,m)=1$ for all type 2 states.
     \end{enumerate}
\end{theorem}

Theorem \ref{theorem:opt_policy_struct} characterizes the structure of the optimal dispatching policy $(\lambda^\ast, \phi^\ast)$ based on the relationship between the driver and rider penalty coefficients $c_d$ and $c_r$. %
Notably, when \( c_d = c_r \)---that is, when the difference between the drivers' in-service and idle costs equals the difference between the rider's queueing and pickup waiting penalties (\( w_s^d - w_o^d = w_q^r - w_p^r \))---the optimal dispatching policy admits a closed-form solution: after each dispatch, the system should be positioned along the boundary between type 1 and type 2 states, at the state with the highest trip completion rate on the corresponding anti-diagonal.
When \( c_d \neq c_r \), the optimization simplifies to determining the dispatching policy for only one state type. Specifically, when \( c_d < c_r \), the rider queueing penalty has a greater impact, so the platform should adopt an aggressive stance by dispatching drivers more readily to reduce rider wait times. Conversely, when \( c_d > c_r \), driver penalty dominates, and the platform should take a more conservative approach by delaying dispatches and letting more riders queue in order to reduce pick-up distances and associated driving costs. The proof relies on a series of coupling arguments, where we couple the system under an arbitrary policy with one governed by a deliberately constructed policy that satisfies the properties in Theorem \ref{theorem:opt_policy_struct}. We then show that the constructed policy leads to a weakly better objective. Although one state type is fully characterized when $c_d \neq c_r$, it is challenging to derive an optimal policy for the remaining state type. When $c_d>c_r$, within type 2 states, dispatching a driver will increase the trip completion rate (and hence the future revenue). However, it also moves the system from $(l,m)$ to $(l+1,m-1)$, so the penalty rate changes from $c_dl+c_rm$ to $c_d(l+1)+c_r(m-1)$. Because $c_d(l+1)+c_r(m-1)-(c_dl+c_rm)=c_d-c_r>0$, this dispatch decision raises the penalty rate. The platform therefore faces an important trade-off: it is not immediately clear whether the additional revenue generated by a faster trip completion rate outweighs the higher penalties induced by dispatching. These effects are not directly comparable, because the benefit of faster trip completion is realized in the future, whereas the penalty rate increases immediately. The same issue arises in type 1 states when $c_d<c_r$, where dispatching will reduce penalty costs but at the expense of lower trip completion rate. Because the revenue and penalty effects of dispatching move in opposite directions and depend on the state, the optimal dispatching policy does not admit a simple cutoff structure across type 1 and type 2 states.

We conclude this subsection by noting that Theorem~\ref{theorem:opt_policy_struct} also holds under any static pricing policy, where the effective request rate remains constant across all states.

\begin{remark}{\label{remark:theorem_1_static_price}} 
If the effective request rate is held constant across all states (i.e., $\lambda(l,m) = \bar{\lambda}$ for all $l \in \mathcal{L}$ and $m \in \mathbb{Z}_{\ge 0}$), an optimal dispatching policy still satisfies Theorem~\ref{theorem:opt_policy_struct}.
\end{remark}

\subsection{Zigzag Dispatching Policy}\label{sec:zigzag_policy} %

When $c_d = c_r$, the optimal dispatching policy $\phi^*$ described in Theorem~\ref{theorem:opt_policy_struct} is characterized by a zigzag boundary that separates the regions where $\phi^*(l, m) = 0$ and $\phi^*(l, m) = 1$ (see Table~\ref{subtab:zigzag_policy_example} for an example). This observation motivates us to define a specific subclass of dispatching policies, which we refer to as \emph{zigzag} policies. As we will show in subsequent sections, even without the condition $c_d = c_r$ or Assumption~\ref{as:2}, focusing on this subclass can substantially simplify computation while still yielding near-optimal solutions. We begin by formally defining a zigzag policy.

\smallskip

\begin{definition}[Zigzag Policy]\label{def:zigzag}
We say a dispatching policy is zigzag if:
\begin{enumerate}
    \item[(1)] For all $l\in\mathcal{L}$, $(\phi(l,m))_{m\in\mathbb Z_{\ge0}}$ is a sequence of consecutive $0$(s) followed by consecutive $1$(s);
    \item[(2)] For all $m\in\mathbb Z_{\ge0}$, $(\phi(l,m))_{l\in\mathcal{L}}$ is a sequence of consecutive $1$(s) followed by consecutive $0$(s).
\end{enumerate}
\end{definition}

\smallskip

If we interpret the policy as a matrix with rows indexed by $l$ and columns by $m$, the definition of a zigzag policy closely resembles that of a matrix in row echelon form: all elements to the left of the leading entry in each row are zero, and each leading entry appears to the right of the leading entry in the row above. The key difference is that, in a zigzag policy, all nonzero entries are set to 1. Below, we provide examples of both zigzag and non-zigzag policies for illustration.

\setlength{\tabcolsep}{6pt} %
\renewcommand{\arraystretch}{1.5}
\begin{table}[htbp]
    \small
    \begin{subtable}{.33\linewidth}
      \centering

        \begin{tabular}{c|cccccc}
         \backslashbox[10mm]{$l$}{$m$} & 0 & 1 & 2 & 3 & 4 & $\cdots$ \\
         \hline 0 & 0 & \textbf{0}\tikzmark{p1} & 1 & 1 & 1 & $\cdots$ \\
         1 & 0\tikzmark{p2} & \textbf{0}\tikzmark{p3} & 1 & 1 & 1 & $\cdots$ \\
         2 & 0 & \textbf{0}\tikzmark{p4} & \textbf{0}\tikzmark{p5} & \textbf{0}\tikzmark{p6} & 1 & $\cdots$ \\
         3 & 0 & 0 & 0 & \textbf{0}\tikzmark{p7} & \textbf{0}\tikzmark{p8} & $\cdots$ \\
         4 & 0 & 0 & 0 & 0 & 0 & $\cdots$ \\
        \end{tabular}
        \begin{tikzpicture}
        [overlay, remember picture, transform canvas={yshift=.05\baselineskip,xshift=-.25\baselineskip},color=black]
        
        \draw [thick,-stealth,shorten >= 5pt,color=gray] ({pic cs:p1}) [left] to ({pic cs:p3});

        \draw [thick,-stealth,shorten >= 5pt,color=gray] ({pic cs:p3}) [left] to ({pic cs:p4});

        \draw [thick,-stealth,shorten >= 5pt,color=gray] ({pic cs:p6}) [left] to ({pic cs:p7});

        \end{tikzpicture}

        \begin{tikzpicture}
        [overlay, remember picture, transform canvas={yshift=.25\baselineskip,xshift=-0.1\baselineskip},color=black]
        \draw [thick,-stealth,shorten >= 3pt,color=gray] ({pic cs:p4}) [left] to ({pic cs:p5});
        \draw [thick,-stealth,shorten >= 3pt,color=gray] ({pic cs:p5}) [left] to ({pic cs:p6});
        \draw [thick,-stealth,shorten >= 3pt,color=gray] ({pic cs:p7}) [left] to ({pic cs:p8});
        \end{tikzpicture}
        \caption{Zigzag policy}
        \label{subtab:zigzag_policy_example}
    \end{subtable}%
    \begin{subtable}{.33\linewidth}
      \centering
        
         \begin{tabular}{c|cccccc}
         \backslashbox[10mm]{$l$}{$m$} & 0 & 1 & 2 & 3 & 4 & $\cdots$ \\
         \hline 0 & 0 & 1 & 1 & 1 & 1 & $\cdots$ \\
         1 & 0 & 0 & 1 & 0 & 1 & $\cdots$ \\
         2 & 0 & 0 & 0 & 0 & 1 & $\cdots$ \\
         3 & 0 & 0 & 0 & 0 & 0 & $\cdots$ \\
         4 & 0 & 0 & 0 & 0 & 0 & $\cdots$ \\
        \end{tabular}
        \caption{Non-zigzag policy}
        \label{subtab:nonzigzag_1}
    \end{subtable} 
    \begin{subtable}{.33\linewidth}
      \centering
        
         \begin{tabular}{c|cccccc}
         \backslashbox[10mm]{$l$}{$m$} & 0 & 1 & 2 & 3 & 4 & $\cdots$\\
         \hline 0 & 0 & 0 & 1 & 1 & 1 & $\cdots$\\
         1 & 0 & 1 & 1 & 1 & 1 & $\cdots$ \\
         2 & 0 & 0 & 0 & 0 & 1 & $\cdots$\\
         3 & 0 & 0 & 0 & 0 & 0 & $\cdots$\\
         4 & 0 & 0 & 0 & 0 & 0 & $\cdots$ \\
        \end{tabular}
        \caption{Non-zigzag policy}\label{subtab:nonzigzag_2}
    \end{subtable} 
    
    \caption{Examples of zigzag and non-zigzag policies.}
    \label{tab:zigzag_example}
\end{table}

\begin{example}[Zigzag and non-zigzag policies]
    In Table \ref{tab:zigzag_example}, we consider a platform with a total of $L=4$ drivers. For $m>4$, we keep the dispatching policy the same as the column with $m = 4$. Each element represents whether to make a dispatch under the corresponding state. Table \ref{subtab:zigzag_policy_example} gives an example of a zigzag policy that meets both conditions outlined in Definition \ref{def:zigzag}. In contrast, the policies depicted in Tables \ref{subtab:nonzigzag_1} and \ref{subtab:nonzigzag_2} are not zigzag policies. The policy depicted in Table \ref{subtab:nonzigzag_1} violates condition (1), while the policy depicted in Table \ref{subtab:nonzigzag_2} violates condition (2) in the definition.
\end{example}

\smallskip

Note that in the absence of Assumption~\ref{as:2} and $c_r=c_d$, 
there might not exist a zigzag policy within the set of optimal dispatching policies. In Appendix \ref{subsec:add_ex}, we provide a counterexample (Example \ref{ex:sub_opt_zigzag}). %
Nevertheless, focusing on zigzag policies offers both practical and analytical advantages. First, zigzag policies are easily interpretable as threshold policies: for a given number of drivers in service $l$, the policy specifies a dispatching threshold on the rider queue length $m$, so the platform dispatches whenever $m$ exceeds this threshold. Conversely, for a given $m$, the policy sets a threshold on $l$, dispatching whenever $l$ falls below this value. Second, zigzag policies are highly tractable, as the resulting Markov chain reduces to a birth-death process that enables closed-form evaluation of the objective. In our model, a birth refers to a rider joining the queue, and a death refers to a trip completion. Moreover, under a zigzag policy, the set of recurrent states traces out a \emph{zigzag path} starting at a state with $l=0$. This path serves as a boundary between regions of zeros and ones, for example, $((0,1), (1,1), (2,1), (2,2), (2,3), (3,3),(3,4),\cdots)$ as highlighted in Table \ref{subtab:zigzag_policy_example}. Since only recurrent states affect the system's long-run behavior, we can focus our analysis only on states along a zigzag path. We now formally introduce the definition of a zigzag path, which will play a key role in the algorithm design. %

\smallskip

\begin{definition}[Zigzag Path]\label{def:path} A zigzag path $P=((l_1,m_1),(l_2,m_2),\cdots,(l_I,m_I))$ under a zigzag dispatching policy $\phi$ is a sequence of states that satisfies:%
\begin{enumerate}
    \item[(1)] (Origin) $l_1=0$, $\phi(0,m_1)=0$ and $\phi(0,m_1+1)=1$;
    \item[(2)] (Monotonicity) $l_{i+1}\geq l_i, m_{i+1}\geq m_i$, $\forall i\in\{1,\cdots, I-1\}$;%
    \item[(3)] (Single movement) $l_{i+1}+m_{i+1}-l_i-m_i=1$, $\forall i\in\{1,\cdots, I-1\}$.
\end{enumerate}
\end{definition}

\smallskip

In other words, a path is defined as a sequence of states, where each successive state is either one step to the right (i.e., $l_{i+1}=l_i$ and $m_{i+1}=m_i+1$), or one step down (i.e., $l_{i+1}=l_i+1$ and $m_{i+1}=m_i$) from the previous state. Given index $i$, the rightward step indicates $\phi(l_i,m_i+1)=0$, which means that the platform will not make a dispatch for a new rider arrival. However, it will dispatch a driver for a trip completion. Conversely, the downward step indicates $\phi(l_i,m_i+1)=1$ and $\phi(l_i+1,m_i)=0$, which means that the platform will not make a dispatch upon a trip completion but will dispatch a driver upon a new rider arrival. %

Each zigzag path $P = ((l_1, m_1), (l_2, m_2), \ldots, (l_I, m_I))$ defines a zigzag policy that terminates at state $(l_I, m_I)$ by rejecting any new riders. The resulting stationary distribution can be computed in closed form using standard results for birth-death processes. %
\begin{equation}\label{eq:terminal_state_stationary_prob}
\pi(l_i,m_i)=
\begin{cases}
\displaystyle
\left(\sum_{k=1}^{I}\prod_{j=1}^{k-1}
\frac{\lambda(l_j,m_j)}{l_{j+1}\mu_{l_{j+1},m_{j+1}}}\right)^{-1},
& i=1,\\[1em]
\displaystyle
\frac{\prod_{j=1}^{i-1}
\frac{\lambda(l_j,m_j)}{l_{j+1}\mu_{l_{j+1},m_{j+1}}}}
{\sum_{k=1}^{I}\prod_{j=1}^{k-1}
\frac{\lambda(l_j,m_j)}{l_{j+1}\mu_{l_{j+1},m_{j+1}}}},
& i=2,\ldots,I.
\end{cases}
\end{equation}

\smallskip

{\textbf{A Dynamic Programming (DP) Heuristic for Optimal Zigzag Policy.}} Leveraging the properties of zigzag paths, we now present an efficient dynamic-programming-based heuristic algorithm that produces high-quality zigzag dispatching policies for general settings with arbitrary $c_d, c_r > 0$, without requiring Assumption~\ref{as:2}.
The algorithm starts at $(0,0)$, iterates top to bottom through the first column ($m=0$), and then proceeds column by column, iteratively extending the zigzag path. %
When reaching an unvisited state $(l, m)$, the algorithm evaluates all zigzag paths terminating at that state, considering transitions from two possible directions: (1) from the left, where the previous state is $(l, m-1)$ (if $m \geq 1$), and (2) from above, where the previous state is $(l-1, m)$ (if $l \geq 1$). Enumerating all possible paths results in an exponential number of candidates in $l$ and $m$. To address this, we introduce a \emph{dominance condition} to prune suboptimal paths at each state. Specifically, in our framework, a zigzag path is said to dominate another terminating at the same state if, under the corresponding optimal pricing policy, its objective value is always greater for any future extensions. Consequently, dominated paths can be safely discarded, as they cannot contribute to the optimal zigzag dispatching policy. For example, a sufficient set of exact dominance conditions for two zigzag paths terminating at the same state, path 1 (weakly) dominates path 2, can be:
\begin{itemize}
    \item[(1)] \emph{Service rate dominance.} For every state $(l,m)$ in path 2, the trip completion rate $l\mu_{l,m}$ must be no greater than that of the state on path 1 lying on the same anti-diagonal (i.e., the state with the same $l+m$ value);
    \item[(2)] \emph{Penalty rate dominance.} For every state $(l,m)$ in path 2, the total penalty rate $c_dl + c_r m$ must be greater than or equal to that of the state on path 1 lying on the same anti-diagonal (i.e., the state with the same $l+m$ value).
\end{itemize}
A formal statement of this dominance condition is given in Proposition~\ref{prop:path_comparison} in Appendix~\ref{sec:perf_evaluation_alg_1}. Although these criteria guarantee dominance, they are too stringent to be satisfied in numerical experiments.
Meanwhile, it is challenging to establish an exact dominance condition that guarantees only a polynomial or constant number of paths are retained at each state. 
This difficulty arises because our objective is inherently path-dependent: when a zigzag path is extended by additional states, the marginal change in its objective value (under the optimal pricing policy) depends on the induced stationary distributions, which in turn depend on the entire sequence of states traversed (see equation \eqref{eq:terminal_state_stationary_prob}), rather than just the newly added state. This is in sharp contrast to stage-additive shortest-path problems, where edge costs are simply summed and classical dominance rules readily apply (such as those in Dijkstra's algorithm); such rules do not hold in our setting.

Instead, we adopt a \emph{heuristic} dominance rule: if one zigzag path achieves a higher objective value than another terminating at the same state, it is likely, though not guaranteed, that this dominance will persist when both paths are extended with identical future states.\footnote{In Appendix~\ref{sec:perf_evaluation_alg_1}, we report an empirical validation to assess this heuristic. To numerically quantify the likelihood that such dominance persists, we measure how often the policies using such heuristic dominance rule are exactly optimal. We find this frequency is relatively high ($\ge$ 60\%) for $L<10$, suggesting that the heuristic works well in small instances.} Based on this heuristic, at each state, we retain only the path with the higher objective value among those arriving from the left or from above, without considering the influence of future states.
On the other hand, when we compute the objective of each zigzag path, its optimal pricing policy can be found by solving the corresponding Bellman equation using value iteration.\footnote{One could also directly optimize it using the closed-form expressions from equation \eqref{eq:terminal_state_stationary_prob} through a nonlinear optimization solver.} Although its stationary distribution admits a closed-form expression from equation \eqref{eq:terminal_state_stationary_prob}, this approach can be computationally intensive, especially when used as a subroutine, since the prices for every state $(l, m)$ along the zigzag path $P$ must be optimized.
To maintain efficiency, \emph{static pricing} is often employed as an effective approximation (see, e.g., \citealt{besbes2022static} and \citealt{bergquist2024staticpricing}). Let $\tilde{\mathcal{R}}_s(\bar{\lambda},P):=\sum_{(l,m)\in P} \pi(l,m)\left(\bar{\lambda} \left(p_0 + w_p^rt_0 + p_1(\bar{\lambda})\cdot t_0\right)- c_d l - c_r m\right)$ be the objective value under zigzag path $P$ and static arrival rate $\bar{\lambda}$. This quantity can be computed and optimized (using line search over $\bar{\lambda}$) efficiently using the stationary distribution $\pi(l,m)$ from equation \eqref{eq:terminal_state_stationary_prob}. %

The entire procedure is summarized in Algorithm~\ref{alg:DP_zigzag}. We construct a table $P$, where $P_{l,m}$ records the best zigzag path terminating at state $(l,m)$. Notably, the update steps in Algorithm~\ref{alg:DP_zigzag} naturally incorporate the option to terminate a zigzag path early by blocking all future arrivals beyond a certain queue length---implemented by setting the maximum price $p_{\mathrm{max}}$.\footnote{Extending a path may decrease the objective value under static pricing, as the constraint of uniform pricing across all states can lead to large penalty terms. In such cases, the optimal static pricing may be a poor approximation of the optimal dynamic policy. Following the approach of \cite{bergquist2024staticpricing}, we consider applying the static price $\bar{\lambda}$ only up to a designated cutoff state along the zigzag path. Beyond this cutoff, all subsequent arrivals are blocked by setting the maximum price $p_1=p_{\mathrm{max}}$.} Accordingly, we maintain a matrix $R$, where $R_{l,m}$ stores the best objective value under static pricing for zigzag paths that either terminate at $(l,m)$ or are cut off prior to reaching $(l,m)$. 
At each state $(l,m)$, let $P_{\text{above}}$ and $P_{\text{left}}$ denote the zigzag paths that reach the current state by extending the best paths from the state above and the state to the left, respectively.
Specifically, consider extending the path $P_{l-1,m}$ vertically downward by one state to $(l,m)$, denoted as $P_{\text{above}}$. Let $\mathcal{R}_{\text{above}}$ be the maximum of the optimal objective value achieved by $P_{\text{above}}$ under static pricing and $R_{l-1,m}$. If $R_{l-1,m}$ is larger, this indicates that introducing a cutoff yields a higher objective, so extending $P_{l-1,m}$ to $P_{\text{above}}$ does not improve the objective value, and thus $\mathcal{R}_{\text{above}} = R_{l-1,m}$. Conversely, if $R_{l-1,m}$ is smaller, no cutoff is imposed, and optimizing the static price along $P_{\text{above}}$ increases the objective value; in this case,
$\mathcal{R}_{\text{above}} = \sup_{\bar{\lambda}} \tilde{\mathcal{R}}_s(\bar{\lambda}, P_{\text{above}})$. The same reasoning applies for $P_{\text{left}}$. After that, we retain the path between $P_{\text{above}}$ and $P_{\text{left}}$ that has a higher objective value for $P_{l,m}$. If two paths have the same objective value, we break ties by selecting the one with the higher trip completion rate---that is, we choose $P_{\text{above}}$ if the trip completion rate at $(l-1, m)$ exceeds that at $(l, m-1)$ (i.e., if $(l-1, m)$ is a type 1 state); otherwise, we retain $P_{\text{left}}$.

As we briefly mentioned in Section \ref{sec:model_and_prelim}, although the number of riders waiting in the queue can grow without limit, %
there naturally exists an upper bound $M$ without loss of optimality due to the increasing rider waiting penalty. Proposition \ref{prop:M_limit} below gives such a bound in the context of Algorithm \ref{alg:DP_zigzag} to facilitate an implementation with a finite number of iterations $\mathcal{O}(ML)$. 

\begin{proposition}
\label{prop:M_limit}
If $M \geq L \bar{\mu} (p_0+p_{\mathrm{max}}t_0) / c_r$, then Algorithm~\ref{alg:DP_zigzag} always returns policies with the same objective value.
\end{proposition}

After the algorithm terminates, we compute the optimal dynamic pricing policy $\lambda^\ast$ using value iteration based on the dispatching policy generated by the zigzag path $P_{L,M}$. Theorem~\ref{theorem:dp_struct} demonstrates that, although Algorithm~\ref{alg:DP_zigzag} is heuristic in nature, it recovers the optimal dispatching and pricing policies under Assumption~\ref{as:2} and $c_d = c_r$.

\begin{algorithm}[htbp] \small \caption{A DP Heuristic for Optimal Zigzag Policy}\label{alg:DP_zigzag} 
\begin{algorithmic}[1] 
\small
\State Create an $(L+1) \times (M+1)$ table $P$, where $P_{l,m}$ stores the best zigzag path terminated at $(l,m)$.
\State Create an $(L+1) \times (M+1)$ matrix $R$, where $R_{l,m}$ stores the best objective value at $(l,m)$.
\State Initialize $P_{l,0}\leftarrow ((0,0),(1,0),\cdots,(l,0))$ for all $l\in\{0,1,\cdots,L\}$; $P_{0,m}\leftarrow ((0,m))$ for all $m\in\{1,2,\cdots,M\}$.
\State Initialize $R_{0,m} \leftarrow 0$ for all $m\in\{0,1,\cdots, M\}$.
\For{$l=1$ to $L$} 
\State $R_{l,0}\leftarrow \max\{R_{l-1,0},~  \sup_{\bar{\lambda}} \tilde{\mathcal{R}}_s(\bar{\lambda},P_{l,0})\}$.
\EndFor

\For{$m=1$ to $M$} \For{$l=1$ to $L$} 
\State Construct two paths $P_{\text{above}}$ and $P_{\text{left}}$ by appending $(l,m)$ to $P_{l-1,m}$ and $P_{l,m-1}$, respectively.
\State Compute $\mathcal{R}_{\text{above}} \leftarrow 
\max\{R_{l-1,m},~\sup_{\bar{\lambda}} \tilde{\mathcal{R}}_s(\bar{\lambda},P_{\text{above}})\}$.
\State Compute $\mathcal{R}_{\text{left}} \leftarrow
\max\{R_{l,m-1},~\sup_{\bar{\lambda}} \tilde{\mathcal{R}}_s(\bar{\lambda},P_{\text{left}})\}$. \If{$\mathcal{R}_{\text{left}}> \mathcal{R}_{\text{above}}$} 
\State  $R_{l,m}\leftarrow\mathcal{R}_{\text{left}}$, $P_{l,m} \leftarrow P_{\text{left}}$.
\ElsIf {$\mathcal{R}_{\text{left}}< \mathcal{R}_{\text{above}}$} 
\State  $R_{l,m}\leftarrow\mathcal{R}_{\text{above}}$, $P_{l,m} \leftarrow P_{\text{above}}$.
\ElsIf {$\mathcal{R}_{\text{left}}= \mathcal{R}_{\text{above}}$}
\State $R_{l,m}\leftarrow\mathcal{R}_{\text{left}}$, 
$P_{l,m} \gets
\begin{cases}
  P_{\text{above}}, & \text{if } (l-1,m)\text{ is a type 1 state}.\\
  P_{\text{left}},  & \text{if } (l-1,m)\text{ is a type 2 state}.
\end{cases}$
\EndIf \EndFor \EndFor 
\State Compute the optimal dynamic pricing policy $\lambda^\ast$ under zigzag path $P_{L,M}$.
\State \textbf{Output:} The dispatching policy yielded by zigzag path $P_{L,M}$ and the corresponding pricing policy $\lambda^\ast$. 
\end{algorithmic} 
\end{algorithm}

\begin{theorem}\label{theorem:dp_struct}
Suppose $\mu$ satisfies Assumption~\ref{as:2}, and $c_d = c_r$. Algorithm~\ref{alg:DP_zigzag} yields optimal dispatching and pricing policies $(\phi^*, \lambda^*)$ such that $\phi^*(l, m) = 0$ for all type 1 states and $\phi^*(l, m) = 1$ for all type 2 states.
\end{theorem}

Finally, in Appendix~\ref{subsec:pricing}, we investigate structural properties of the optimal pricing policy. In particular, we characterize conditions under which the optimal price is nondecreasing along any zigzag path. This monotone pricing structure resembles the \emph{surge pricing} commonly observed on ride-hailing platforms: when the number of nearby idle drivers is low, the platform increases rider prices to better balance supply and demand. We extend this intuition to a two-sided spatial queueing system in which riders may also wait in queue.

\section{Numerical Experiments}\label{sec:numerical_experiments}

In this section, we demonstrate the effectiveness and efficiency of our proposed Algorithm \ref{alg:DP_zigzag} through extensive numerical experiments. In Section \ref{subsec:synthetic_exp}, we conduct synthetic experiments within our continuous-time Markovian model. %
In Section \ref{subsec:simulation}, we evaluate a range of dispatching and pricing policies in a simulated environment that incorporates real-world spatial features, relinquishing various modeling assumptions.

\subsection{Synthetic Experiments}\label{subsec:synthetic_exp}
We conduct experiments based on our Markovian model, considering two scenarios with parameters $(L, \Lambda, M) \in \{(20, 8, 10), (100, 40, 50)\}$. Rider arrival locations and destinations are independently and uniformly sampled from a $10 \times 10$ square service region, with distances calculated using the Euclidean norm. All distances and times are measured in kilometers and minutes, respectively. We assume drivers travel in unit speed. The service rate $\mu_{l,m}$ is estimated via Monte Carlo sampling: specifically, we uniformly sample the locations of $L-l$ idle drivers and $m$ riders within the service region, and, at each dispatch event, employ a nearest-neighbor strategy that matches the closest rider-driver pair. This sampling procedure is repeated 100,000 times for each $(l, m)$, and the average pickup time is computed.
Finally, we estimate $\mu_{l,m}$ as the reciprocal of the sum of the average trip time and the average pickup time.\footnote{Our estimated service rates satisfy Assumption~\ref{as:2} for all states $(l,m)$ except $(1,1)$ and $(0,2)$. Since these states have negligible stationary probabilities, their exceptions have little impact on optimization.} The average trip time $t_0$ is taken as the expected distance between two uniformly sampled points in the $10 \times 10$ square, which can be computed as $t_0= \frac{2}{3}\left(2+\sqrt{2}+5\ln(1+\sqrt{2})\right)\approx 5.2$. We specify the inverse demand function using a linear model, $p_1(\lambda) = 2(1 - \lambda/\Lambda)$, and set $p_0 + w_p^r t_0 = 5$. In this setting, the base fare per trip $p_0$ is approximately \$4, the per-minute penalty for rider pickup waiting is $w_p^r = 0.2$, and the maximum per-kilometer price $p_{\text{max}}$ is \$2---consistent with those commonly observed on ride-hailing platforms.

We test $3\times3=9$ combinations of $c_d$ and $c_r$ by taking $ c_d\in \{0.5, 0.75, 1.0\}$ and $c_r \in \{0.5, 0.75, 1.0\} $. In each scenario, we evaluate the performance of our zigzag policy, optimized using Algorithm \ref{alg:DP_zigzag}, by comparing it against two benchmark policies.

\begin{enumerate}
    \item \textbf{Value Iteration}: We apply value iteration to solve objective \eqref{eq:objective_alt} using the alternative action space introduced at the beginning of Section~\ref{sec:methodology}. As previously discussed, the full action space, which permits dispatching multiple drivers simultaneously, is computationally intractable. Therefore, we restrict the action space to allow dispatching at most one driver at a time. %
    More details can be found in Appendix \ref{subsec:alt_def_action_space}. 
    \item \textbf{Greedy Dispatching Policy}: A baseline policy that always dispatches a driver whenever one is available, with the dynamic pricing policy subsequently optimized via value iteration.
\end{enumerate}

\begin{table}[htbp]
\SingleSpacedXI
\footnotesize
    \centering
    \begin{tabular}{ccccccccc}
        \toprule 
         & \multicolumn{2}{c}{Value Iteration} & \multicolumn{2}{c}{Greedy, Dyn.} & \multicolumn{2}{c}{Zigzag, Dyn.} & \multicolumn{2}{c}{Zigzag, Sta.}\\
        \cmidrule(lr){2-3} \cmidrule(lr){4-5} \cmidrule(lr){6-7} \cmidrule(lr){8-9}
        $c_d, c_r$ & $\Tilde{\mathcal{R}}(\lambda,\phi)$ & Runtime (s) & $\Tilde{\mathcal{R}}(\lambda,\phi)$ & Runtime (s) & $\Tilde{\mathcal{R}}(\lambda,\phi)$ & Runtime (s) & $\Tilde{\mathcal{R}}_s(\lambda,\phi)$ & Runtime (s) \\
        \midrule
        0.5, 0.5  & 20.54 & 544.79 & 18.32 & 20.39 & 20.54 & 16.56 & 20.11 & 1.69 \\
        0.5, 0.75 & 19.67 & 514.68 & 17.53 & 26.11 & 19.67 & 15.90 & 19.21 & 1.67 \\
        0.5, 1.0  & 19.05 & 469.02 & 17.25 & 19.03 & 19.05 & 14.73 & 18.53 & 1.52 \\
        \midrule
        0.75, 0.5 & 16.36 & 572.13 & 13.61 & 37.42 & 16.36 & 17.51 & 16.04 & 1.65 \\
        0.75, 0.75 & 15.57 & 525.36 & 13.46 & 23.30 & 15.57 & 16.43 & 15.24 & 1.60 \\
        0.75, 1.0 & 15.02 & 466.05 & 13.45 & 15.77 & 15.02 & 15.08 & 14.65 & 1.46 \\
        \midrule
        1.0, 0.5 & 12.44 & 591.52 & 10.06 & 21.21 & 12.43 & 18.21 & 12.23 & 1.45 \\
        1.0, 0.75 & 11.74 & 511.15 & 10.06 & 14.09 & 11.74 & 16.47 & 11.53 & 1.42 \\
        1.0, 1.0 & 11.26 & 451.86 & 10.06 & 12.32 & 11.26 & 14.79 & 11.03 & 1.27 \\
        \bottomrule
    \end{tabular}
    \caption{Comparison of different algorithms under parameters $(L,\Lambda,M)=(20,8,10)$.}
    \label{tab:results_L20}
\end{table}

\begin{table}[htbp]
\SingleSpacedXI
\footnotesize
    \centering
    \begin{tabular}{ccccccccc}
        \toprule 
         & \multicolumn{2}{c}{Value Iteration} & \multicolumn{2}{c}{Greedy, Dyn.} & \multicolumn{2}{c}{Zigzag, Dyn.} & \multicolumn{2}{c}{Zigzag, Sta.}\\
        \cmidrule(lr){2-3} \cmidrule(lr){4-5} \cmidrule(lr){6-7} \cmidrule(lr){8-9}
        $c_d, c_r$ & $\Tilde{\mathcal{R}}(\lambda,\phi)$ & Runtime (s) & $\Tilde{\mathcal{R}}(\lambda,\phi)$ & Runtime (s) & $\Tilde{\mathcal{R}}(\lambda,\phi)$ & Runtime (s) & $\Tilde{\mathcal{R}}_s(\lambda,\phi)$ & Runtime (s) \\
        \midrule
        0.5, 0.5  & 119.07 & 3,600.00 & 116.10 & 411.35 & 128.17 & 488.77 & 126.96 & 167.91 \\
        0.5, 0.75 & 111.87 & 3,600.00 & 111.20 & 2,780.35 & 125.78 & 658.32 & 124.36 & 168.83 \\
        0.5, 1.0  &  106.99 & 3,600.00 & 111.18 & 485.99 & 124.00 & 445.39 & 122.42 & 151.18 \\
        \midrule
        0.75, 0.5 &  95.41 & 3,600.00 &  91.11 & 1,773.91 & 105.41 & 573.68 & 104.67 & 169.17 \\
        0.75, 0.75 &  88.20 & 3,600.00 &  90.80 & 412.81 & 103.33 & 525.72 & 102.46 & 168.77 \\
        0.75, 1.0  &  83.55 & 3,600.00 &  90.80 & 266.46 & 101.76 & 460.57 & 100.80 & 152.40 \\
        \midrule
        1.0, 0.5  &  72.03 & 3,600.00 &  71.37 & 461.67 &  83.52 & 554.08 &  83.18 & 150.50 \\
        1.0, 0.75 &  64.85 & 3,600.00 &  71.37 & 273.75 &  81.76 & 527.70 &  81.35 & 150.32 \\
        1.0, 1.0  &  59.58 & 3,600.00 &  71.37 & 283.30 &  80.43 & 475.77 &  79.97 & 130.77 \\
        \bottomrule
    \end{tabular}
    \caption{Comparison of different algorithms under parameters $(L,\Lambda,M)=(100,40,50)$. \normalfont{A 60-minute cap is imposed on every run; any runtime that exceeds this limit is reported as 3600.00 seconds.}}
    \label{tab:results_L100}
\end{table}

Tables~\ref{tab:results_L20} and~\ref{tab:results_L100} present results for two scenarios, $(L, \Lambda, M) \in \{(20, 8, 10), (100, 40, 50)\}$, respectively. We evaluate two pricing strategies derived from the zigzag policy produced by Algorithm~\ref{alg:DP_zigzag}: a static pricing policy (``Zigzag, Sta.'') and a dynamic pricing policy (``Zigzag, Dyn.''). Both zigzag policies consistently outperform the greedy benchmark. In Table~\ref{tab:results_L20}, where value iteration can reach optimality within the time limit, ``Zigzag, Dyn.'' attains nearly identical objective values but with substantially lower computation time. Notably, when $c_d = c_r$, ``Zigzag, Dyn.'' matches the objective value of value iteration exactly, confirming its optimality as established in Theorem~\ref{theorem:dp_struct}. For the larger scenario in Table~\ref{tab:results_L100}, value iteration fails to converge within the time limit and yields solutions inferior to even the greedy policy. In Appendix~\ref{sec:add_perf_metric_sync}, we present and discuss additional performance metrics, including the average price per rider, revenue rate, average pickup time, and average queueing time, under the same experimental setup. In Appendix~\ref{sec:perf_evaluation_alg_1}, we conduct additional experiments to further assess the optimality of Algorithm~\ref{alg:DP_zigzag} by benchmarking it against the optimal joint dispatching and pricing policy. %

\subsection{Simulation Experiments}\label{subsec:simulation}
We further validate our findings by implementing the algorithm in a simulated environment that reflects real-world ride-hailing spatial dispatching. In this simulation, we set $(L, \Lambda, M) = (100, 40, 50)$ and again consider a $10 \times 10$ square region where rider pick-up and drop-off locations are uniformly distributed. Idle drivers remain at their drop-off locations and do not reposition. Departing from the previous synthetic experiments, which assumed independently and state-dependent exponentially distributed service times, this simulation introduces service times that are correlated across dispatches due to driver and waiting rider locations and are endogenously sampled by the actual matching distance, which itself depends on the dispatching and pricing policies.

\subsubsection{Estimation of Service Rate.}\label{subsubsec:estimation_serv_rate}
Instead of estimating the service rate $\mu$  based on randomly generated rider and driver locations, as in Section~\ref{subsec:synthetic_exp}, we now estimate $\mu$ directly from dispatching outcomes observed in simulation.
We begin by simulating the system under a \emph{constant-radius} dispatching policy: only the closest rider-driver pair within a radius $r$ is dispatched, generating samples for analysis. A static pricing policy with a fixed effective request rate $\lambda(l,m) = 12$ is applied uniformly across all states. We conduct 24 simulation runs with dispatching radii $r \in \{0.5, 1, 1.5, \dots, 12\}$, each over a time period of $T = 20,000$. During each run, we record the pick-up times for every dispatch. Specifically, whenever the system enters state $(l,m)$ and the next event is either an arrival or a trip completion immediately followed by a dispatch, we attribute the corresponding pickup time to state $(l,m)$. The average pickup time for each state is then computed across all simulation runs, with each average used as a data point for regression analysis. To mitigate sampling noise, states with fewer than 10 recorded samples are excluded from the dataset.

Following \cite{wang2022ondemand}, we model the average pickup time as a power-law function of the number of waiting riders in the queue and the number of idle drivers. Specifically, the average pickup time $\eta_{l,m}$ in state $(l,m)$ is estimated as:\footnote{There exists an expected pickup time at states where $l = L$ or $m = 0$ and $l>0$, as some riders are still being served in these states. We thus add 1 to both $m$ and $L - l$ in equation \eqref{eq:simulation_mu}.}
\begin{equation}\label{eq:simulation_mu}
    \eta_{l,m} = C\cdot (L-l+1)^{\alpha_1}\cdot (m+1)^{\alpha_2}
\end{equation}
where $C > 0$ and $\alpha_1, \alpha_2 < 0$ are some constants. The right-hand side of Equation~\eqref{eq:simulation_mu} takes the form of a Cobb--Douglas production function~\citep{charles1928theory}.
We then perform a log-log regression. %
In specific, the regression is based on the mean function
$
    \log \eta_{l,m} = \log C + \alpha_1\log(L-l+1)+\alpha_2 \log (m+1)
$.

\begin{table}[htbp]
\SingleSpacedXI
\small
    \centering
    \begin{tabular}{cccccc}
        \toprule
         & coef & std. error  & $\mathrm{P}>|\mathrm{t}|$ & [0.025 & 0.975] \\
        \midrule
        $\hat{C}$  & 3.839 & 0.017  & 0.000 & 3.743 & 3.939\\
        $\hat{\alpha}_1$  & -0.274  & 0.004 & 0.000 & -0.282  & -0.267 \\
        $\hat{\alpha}_2$  & -0.192  & 0.003  & 0.000 & -0.199 & -0.186 \\
        \bottomrule
    \end{tabular}
    \caption{Regression results for parameters $\alpha_1$, $\alpha_2$, and $C$.}
    \label{tab:regression_result}
\end{table}

The regression results are summarized in Table~\ref{tab:regression_result}. We observe that the estimated coefficients, $\hat{\alpha}_1$ and $\hat{\alpha}_2$, are less negative, i.e., smaller in magnitude, than $-0.5$ reported when riders and drivers are uniformly distributed in space (see values reported in \citealt{yan2020dynamic} and \citealt{wang2022ondemand}). This discrepancy is likely attributable to the correlation introduced by driver and waiting rider locations across dispatches. %
The service rate $\mu_{l,m}$ can then be approximated as the reciprocal of the sum of the estimated expected pickup time at state $(l,m)$ and the expected trip time $t_0$, i.e., $
\hat{\mu}_{l,m} = 1/(\hat{\eta}_{l,m} + t_0)$.

\begin{proposition}\label{prop:simulation_mu_property}
    If the fitted coefficient $\hat{C} \leq t_0$, then $\hat{\mu}$ satisfies Assumption \ref{as:2}, and is strictly increasing in $m$ and strictly decreasing in $l$.
\end{proposition}
Proposition~\ref{prop:simulation_mu_property} establishes that $\hat{\mu}$ satisfies our assumptions provided that the expected pickup time at $(l, m) = (0, 0)$ does not exceed $t_0$; that is, the expected pickup time for a single waiting rider and a single available driver must be no greater than the expected trip time $t_0$.
Our regression yields $\hat{C} = 3.839$, which is well below $t_0 \approx 5.2$. Thus, the estimated $\hat{\mu}$ satisfies Assumption~\ref{as:2} with a comfortable margin.

\subsubsection{Simulation Results.}\label{subsubsec:evaluation_of_DP-Heuristic}

We compare the following dispatching and pricing policies over a simulated time period of $T=20,000$:

\begin{enumerate}
    \item \emph{Constant radius dispatching, static pricing}: The platform dispatches the nearest driver-rider pair whenever their distance is less than or equal to a specified threshold $\bar{r}$, and applies a static price (effective request rate) $\bar{\lambda}$. We employ simulation optimization with coordinate ascent to jointly determine the optimal values of $\bar{r}$ and $\bar{\lambda}$.
    \begin{enumerate}
        \item[(a)] We first determine $\bar{r}$ through a line search with a step size of $0.2$ under $\bar{\lambda}=12$ by running the simulation. 
        \item[(b)] Fixing the $\bar{r}$ from step (a), we update $\bar{\lambda}$ using line search with a step size of 0.2 by running the simulation.
        \item[(c)] We repeat the above steps until the absolute differences between consecutive values of $\bar{r}$ and $\bar{\lambda}$ are less than or equal to 0.2.
    \end{enumerate}
    \item \emph{Zigzag dispatching, static pricing}: We adopt the zigzag dispatching policy from Algorithm~\ref{alg:DP_zigzag}, computed using the estimated service rates $\hat{\mu}$ from Section~\ref{subsubsec:estimation_serv_rate}. We use the corresponding optimal static pricing policy based on this dispatching policy.

    \item \emph{Zigzag dispatching, dynamic pricing}: We adopt the same zigzag dispatching policy as in (2), but instead apply the corresponding optimal dynamic pricing policy.

\end{enumerate}

\begin{table}[htbp]
\SingleSpacedXI
\centering
\small
\begin{tabular}{cc c c c c c c}
\toprule
\multicolumn{2}{c}{} 
& \multicolumn{2}{c}{C.R.} 
& \multicolumn{2}{c}{Zigzag, Sta.} 
& \multicolumn{2}{c}{Zigzag, Dyn.} \\
\cmidrule(lr){3-4} \cmidrule(lr){5-6} \cmidrule(lr){7-8}
$c_d$ & $c_r$ 
& Obj. & Price 
& Obj. & Price 
& Obj. & Price \\
\midrule
0.50 & 0.50 & 114.04 & 10.68 & 113.01 & 10.82 & 113.49 & 10.74 \\
0.50 & 0.75 & 111.25 & 10.84 & 110.55 & 10.90 & 111.18 & 10.81 \\
0.50 & 1.00 & 109.45 & 10.94 & 108.89 & 10.97 & 109.68 & 10.86 \\
\midrule
0.75 & 0.50 &  90.91 & 10.84 &  90.67 & 10.90 &  90.96 & 10.82 \\
0.75 & 0.75 &  89.12 & 10.94 &  88.70 & 10.99 &  89.01 & 10.90 \\
0.75 & 1.00 &  87.85 & 11.05 &  87.38 & 11.05 &  87.78 & 10.95 \\
\midrule
1.00 & 0.50 &  69.68 & 11.05 &  69.43 & 11.04 &  69.55 & 10.98 \\
1.00 & 0.75 &  68.24 & 11.15 &  67.92 & 11.13 &  68.12 & 11.06 \\
1.00 & 1.00 &  67.47 & 11.25 &  67.03 & 11.19 &  67.20 & 11.11 \\
\bottomrule
\end{tabular}
\caption{Objective values and average prices under different policies in simulation across various \((c_d,c_r)\) settings.}
\label{tab:combined_objective_pricing}
\end{table}

It is important to note that, in contrast to the constant radius policy, computing both zigzag policies does \emph{not} require access to the simulation environment other than the calibrated service rate $\hat{\mu}$. Table \ref{tab:combined_objective_pricing} shows the objective and average price per ride for various algorithms under different $c_d$ and $c_r$, averaged over the simulation period. In specific,  the average price is computed as $\sum_{(l,m)\in \mathcal{S}} \pi_{\lambda,\phi}(l,m)\left(p_0+p_1(\lambda(l,m))\cdot t_0\right)$. Several observations can be drawn from the table. First, the performance of the two zigzag policies is comparable to that of the constant radius policy, even though the latter directly optimizes its parameters using the same simulation as in the evaluation. Second, dynamic pricing slightly outperforms static pricing under the same zigzag dispatching policy, while also yielding lower average prices. Table~\ref{tab:three_policy_metrics} further reports the average revenue per unit time, average pickup time, and average queueing time under varying values of $c_d$ and $c_r$.

\begin{table}[htbp]
\SingleSpacedXI
\centering
\small
\begin{tabular}{cc ccc ccc ccc}
\toprule
\multicolumn{2}{c}{} 
& \multicolumn{3}{c}{C.R.} 
& \multicolumn{3}{c}{Zigzag, Sta.} 
& \multicolumn{3}{c}{Zigzag, Dyn.}\\
\cmidrule(lr){3-5}\cmidrule(lr){6-8}\cmidrule(lr){9-11}
$c_d$ & $c_r$ 
& Revenue & Pickup & Queue 
& Revenue & Pickup & Queue 
& Revenue & Pickup & Queue \\
\midrule
0.50 & 0.50 & 152.12 & 1.33 & 0.89 & 148.58 & 1.32 & 0.73 & 149.44 & 1.37 & 0.68 \\
0.50 & 0.75 & 147.81 & 1.43 & 0.53 & 145.98 & 1.38 & 0.52 & 147.08 & 1.44 & 0.48 \\
0.50 & 1.00 & 144.89 & 1.50 & 0.36 & 144.17 & 1.42 & 0.41 & 145.50 & 1.50 & 0.37 \\
\midrule
0.75 & 0.50 & 147.92 & 1.35 & 0.60 & 146.09 & 1.27 & 0.63 & 147.32 & 1.31 & 0.60 \\
0.75 & 0.75 & 144.88 & 1.33 & 0.46 & 143.54 & 1.32 & 0.45 & 144.98 & 1.37 & 0.42 \\
0.75 & 1.00 & 141.85 & 1.38 & 0.30 & 141.70 & 1.35 & 0.35 & 143.39 & 1.42 & 0.32 \\
\midrule
1.00 & 0.50 & 141.88 & 1.22 & 0.46 & 141.95 & 1.20 & 0.54 & 143.10 & 1.22 & 0.51 \\
1.00 & 0.75 & 138.67 & 1.27 & 0.30 & 139.38 & 1.24 & 0.39 & 140.82 & 1.27 & 0.36 \\
1.00 & 1.00 & 135.43 & 1.23 & 0.25 & 137.57 & 1.27 & 0.30 & 139.27 & 1.30 & 0.28 \\
\bottomrule
\end{tabular}
\caption{Comparison of different policies in terms of revenue rate, average pickup time, and average queueing time in simulation across various $(c_d, c_r)$ settings.}
\label{tab:three_policy_metrics}
\end{table}

\subsubsection{Robustness Evaluation.}
To further assess the robustness of our algorithm under varying environmental conditions, we evaluate the policies across a range of values for $\Lambda \in \{20, 24, \dots, 60\}$ in the simulation and examine how deviations from the training scenario ($\Lambda=40$) impact policy performance.

\begin{figure}[htbp]
\centering
\begin{tikzpicture}
\begin{groupplot}[
    group style={
        group size=3 by 1,
        horizontal sep=0.35cm
    },
    width=0.295\textwidth,
    height=0.235\textwidth,
    scale only axis,
    grid=both,
    grid style={dotted, gray!50},
    major grid style={dotted, gray!60},
    minor tick num=0,
    tick label style={font=\small},
    label style={font=\small},
    title style={font=\small},
    ylabel style={xshift=2pt},
    legend style={
    font=\tiny,
    draw=gray!25,
    fill=white,
    cells={anchor=west},
    inner xsep=1pt,
    inner ysep=0.5pt,
    row sep=-2pt,
    column sep=1pt
    },
    xlabel={$\Lambda$},
    xmin=18,
    xmax=62,
    ymin=20,
    ymax=140,
    xtick={20,25,30,35,40,45,50,55,60},
    ytick={20,40,60,80,100,120,140}
]

\nextgroupplot[
    title={$(c_1,c_2)=(0.5,0.5)$},
    ylabel={Objective},
    ylabel style={yshift=-7pt},
    legend style={at={(0.97,0.3)}, anchor=north east},
]
\addplot[
    color=mplblue,
    line width=0.8pt,
    solid,
    mark=*,
    mark size=1pt,
    mark options={solid}
] table[col sep=comma, x=Lambda, y=CR] {figures/lambda_obj_panel_050.csv};
\addlegendentry{C.R.}

\addplot[
    color=mplorange,
    line width=0.8pt,
    dashed,
    mark=square*,
    mark size=1pt,
    mark options={solid}
] table[col sep=comma, x=Lambda, y=ZigzagSta] {figures/lambda_obj_panel_050.csv};
\addlegendentry{Zigzag (Sta.)}

\addplot[
    color=mplgreen,
    line width=0.8pt,
    dashdotted,
    mark=triangle*,
    mark size=1pt,
    mark options={solid}
] table[col sep=comma, x=Lambda, y=ZigzagDyn] {figures/lambda_obj_panel_050.csv};
\addlegendentry{Zigzag (Dyn.)}

\nextgroupplot[
    title={$(c_1,c_2)=(0.75,0.75)$},
    yticklabel style={opacity=0},
    legend style={at={(0.97,0.3)}, anchor=north east},
]
\addplot[
    color=mplblue,
    line width=0.8pt,
    solid,
    mark=*,
    mark size=1pt,
] table[col sep=comma, x=Lambda, y=CR] {figures/lambda_obj_panel_075.csv};
\addlegendentry{C.R.}

\addplot[
    color=mplorange,
    line width=0.8pt,
    dashed,
    mark=square*,
    mark size=1pt,
    mark options={solid}
] table[col sep=comma, x=Lambda, y=ZigzagSta] {figures/lambda_obj_panel_075.csv};
\addlegendentry{Zigzag (Sta.)}

\addplot[
    color=mplgreen,
    line width=0.8pt,
    dashdotted,
    mark=triangle*,
    mark size=1pt,
    mark options={solid}
] table[col sep=comma, x=Lambda, y=ZigzagDyn] {figures/lambda_obj_panel_075.csv};
\addlegendentry{Zigzag (Dyn.)}

\nextgroupplot[
    title={$(c_1,c_2)=(1.0,1.0)$},
    yticklabel style={opacity=0},
    legend style={at={(0.97,0.97)}, anchor=north east},
]
\addplot[
    color=mplblue,
    line width=0.8pt,
    solid,
    mark=*,
    mark size=1pt,
] table[col sep=comma, x=Lambda, y=CR] {figures/lambda_obj_panel_100.csv};
\addlegendentry{C.R.}

\addplot[
    color=mplorange,
    line width=0.8pt,
    dashed,
    mark=square*,
    mark size=1pt,
    mark options={solid}
] table[col sep=comma, x=Lambda, y=ZigzagSta] {figures/lambda_obj_panel_100.csv};
\addlegendentry{Zigzag (Sta.)}

\addplot[
    color=mplgreen,
    line width=0.8pt,
    dashdotted,
    mark=triangle*,
    mark size=1pt,
    mark options={solid}
] table[col sep=comma, x=Lambda, y=ZigzagDyn] {figures/lambda_obj_panel_100.csv};
\addlegendentry{Zigzag (Dyn.)}

\end{groupplot}
\end{tikzpicture}
\caption{Robustness of different policies when rider arrival rate $\Lambda$ used in simulation deviates from the training value $\Lambda=40$.}
\label{fig:robust_test}
\end{figure}

Figure~\ref{fig:robust_test} presents the objective values of various algorithms across a range of rider arrival rates in simulation. While all three algorithms exhibit similar performance at the training arrival rate ($\Lambda=40$), their performance diverges as the actual $\Lambda$ deviates from this value. The zigzag policy with dynamic pricing consistently improves its objective value as $\Lambda$ increases and outperforms the other policies across all deviated $\Lambda$ values. In contrast, the other two policies, both relying on static pricing, show greater variability and a decline in objective value when $\Lambda > 40$, despite rising demand. These results underscore the robustness of dynamic pricing, consistent with similar observations in \cite{banerjee2015pricing}. It is also important to note that directly optimizing a dynamic pricing policy $\lambda$ through simulation optimization is intractable due to its high dimension, which further highlights the value of our MDP framework.

\subsubsection{Extension to Multi-region System.}
In practice, a ride-hailing platform often operates across multiple regions with heterogeneous supply and demand conditions. To capture these operational features more faithfully, we extend our simulation environment to a multi-region setting in which drivers may transport riders across regional boundaries and subsequently increase the supply in the destination region.

We model the service area as a $15 \times 15$ square partitioned into a $3 \times 3$ grid of square regions, each of size $5 \times 5$. We label each region by $(i,j)$, where $i,j \in \{1,2,3\}$. To mimic realistic urban mobility patterns, we simulate directional rider flows under two scenarios. In the first scenario, a large fraction of riders originates in peripheral regions, that is, regions with $(i,j) \neq (2,2)$, and travels to the central region $(2,2)$. For each region $(i,j)$, we draw the rider arrival rate $\Lambda_{i,j}$ independently from a $\operatorname{Unif}(4,8)$ distribution. Conditional on the origin region, each rider's destination region is drawn from a categorical distribution with probabilities proportional to region-specific weights: for each peripheral region, the weight is drawn independently from \(\operatorname{Unif}(4,8)\), while the central region is assigned a fixed weight of $10$. In the second scenario, the flow pattern is reversed, with a large fraction of riders traveling from the central region to the eight peripheral regions. Accordingly, we set $\Lambda_{2,2}=10$ for the central region and draw $\Lambda_{i,j}$ independently from a $\operatorname{Unif}(4,8)$ distribution for each peripheral region. Destinations are then sampled uniformly over the service area. In both scenarios, conditional on the region, each location is drawn uniformly within that region.

In such multi-region settings with spatially imbalanced supply and demand, our current implementation may not perform well because it does not account for regional heterogeneity in rider demand and driver availability. To illustrate, consider a two-region system in which most riders travel from region~1 to region~2, while only a small fraction travel in the reverse direction. This directional imbalance causes drivers to accumulate in region~2 and become increasingly scarce in region~1 over time. In this case, the platform should ideally charge a higher price to riders traveling from region~1 to region~2 than to those traveling from region~2 to region~1, since region~1 faces tighter driver supply and longer expected pickup and queueing times due to excess rider demand. Conversely, the platform should lower prices in region~2, where driver supply is abundant and rider waiting times are shorter. More generally, an effective pricing policy in such environments should depend on the regional distribution of both idle drivers and waiting riders. By contrast, the current implementation bases pricing decisions only on the total numbers of idle drivers and riders in the system, and therefore cannot correct for these spatial imbalances.

To address such spatial imbalances in the multi-region setting, we adapt the implementation in Section~\ref{subsec:simulation} by localizing the demand and supply counts. Specifically, rather than using the raw numbers of drivers in service and riders waiting in the queue, we replace these quantities with Gaussian kernel estimates that capture the local density around each arriving rider. When a rider arrives, the platform observes the spatial distribution of idle drivers and assigns each idle driver a weight $K(d)=\exp\!\left(-\frac{1}{2}(d/\sigma)^2\right)$, where $d$ is the Euclidean distance between the arriving rider and the driver, and $\sigma$ is a scale parameter representing the effective neighborhood size. In our implementation, we set $\sigma=5$, equal to the side length of each region. This kernel reflects the intuition that nearby drivers have a greater effect on pickup times and matching feasibility than distant ones. Summing these weights yields a kernel-weighted estimate of local idle supply, and we estimate $l$ as the total number of cars minus this weighted idle count. We apply the same Gaussian kernel construction to waiting riders to obtain a localized measure of demand. Accordingly, when making pricing decisions, the platform uses $\hat{l}:=\max\left\{\,L-l_o\sum_j K(d_j),0\right\}$ and $\hat{m}:=m_o\sum_k K(d_k)$ as spatially smoothed approximations of $l$ and $m$, respectively. In the expression for $\hat{l}$, the index $j$ ranges over all currently idle drivers, and $d_j$ denotes the Euclidean distance between the arriving rider and idle driver $j$, so that $\sum_j K(d_j)$ is a kernel-weighted measure of nearby supply. Similarly, in the expression for $\hat{m}$, the index $k$ ranges over all waiting riders, and $d_k$ denotes the distance from the arriving rider to waiting rider $k$, so that $\sum_k K(d_k)$ provides a localized approximation of nearby demand. The normalization constants $l_o$ and $m_o$ are calibrated offline, using the same simulations employed for service-rate estimation in Section~\ref{subsubsec:estimation_serv_rate}, as the average ratios between the true idle-driver or waiting-rider counts and their corresponding kernel-weighted sums. We round $\hat{l}$ and $\hat{m}$ to the nearest integers. If the resulting state $(\hat{l},\hat{m})$ falls outside the zigzag path, we assign the price associated with the nearest state on the path in the same row.

\begin{table}[htbp]
\SingleSpacedXI
\centering
\small
\begin{tabular}{cc c c c c c c c c}
\toprule
\multicolumn{2}{c}{} 
& \multicolumn{2}{c}{C.R.} 
& \multicolumn{2}{c}{Zigzag, Sta.} 
& \multicolumn{2}{c}{Zigzag, Dyn.} 
& \multicolumn{2}{c}{Zigzag, Dyn. (Loc.)} \\
\cmidrule(lr){3-4} \cmidrule(lr){5-6} \cmidrule(lr){7-8} \cmidrule(lr){9-10}
$c_d$ & $c_r$ 
& Obj. & Price 
& Obj. & Price 
& Obj. & Price 
& Obj. & Price \\
\midrule
0.50 & 0.50 & 172.83 & 14.29 & 172.26 & 14.40 & 172.75 & 14.28 & 174.04 & 14.21 \\
0.50 & 0.75 & 171.01 & 14.29 & 170.11 & 14.52 & 170.83 & 14.37 & 172.15 & 14.29 \\
0.50 & 1.00 & 169.53 & 14.75 & 168.65 & 14.60 & 169.57 & 14.44 & 170.94 & 14.35 \\
\midrule
0.75 & 0.50 & 129.92 & 14.75 & 129.68 & 14.74 & 129.83 & 14.69 & 131.43 & 14.68 \\
0.75 & 0.75 & 129.06 & 14.75 & 128.18 & 14.84 & 128.34 & 14.78 & 130.07 & 14.76 \\
0.75 & 1.00 & 128.37 & 14.75 & 127.14 & 14.92 & 127.32 & 14.84 & 129.08 & 14.82 \\
\midrule
1.00 & 0.50 &  91.38 & 15.21 &  90.92 & 15.37 &  90.93 & 15.36 &  92.79 & 15.38 \\
1.00 & 0.75 &  89.99 & 15.69 &  89.80 & 15.45 &  89.75 & 15.44 &  91.80 & 15.46 \\
1.00 & 1.00 &  89.94 & 15.69 &  89.01 & 15.51 &  88.90 & 15.50 &  91.08 & 15.51 \\
\bottomrule
\end{tabular}
\caption{Objective values and average prices under different pricing and dispatching policies under the demand scenario in which most riders travel from peripheral regions to the central region.}
\label{tab:multi_region_pricing}
\end{table}

\begin{table}[htbp]
\SingleSpacedXI
\centering
\small
\begin{tabular}{cc c c c c c c c c}
\toprule
\multicolumn{2}{c}{} 
& \multicolumn{2}{c}{C.R.} 
& \multicolumn{2}{c}{Zigzag, Sta.} 
& \multicolumn{2}{c}{Zigzag, Dyn.} 
& \multicolumn{2}{c}{Zigzag, Dyn. (Loc.)} \\
\cmidrule(lr){3-4} \cmidrule(lr){5-6} \cmidrule(lr){7-8} \cmidrule(lr){9-10}
$c_d$ & $c_r$ 
& Obj. & Price 
& Obj. & Price 
& Obj. & Price 
& Obj. & Price \\
\midrule
0.50 & 0.50 
& 183.44 & 14.61 
& 183.97 & 14.59 
& 184.47 & 14.50 
& 185.24 & 14.42 \\

0.50 & 0.75 
& 180.95 & 14.61 
& 180.99 & 14.71 
& 181.71 & 14.60 
& 182.70 & 14.50 \\

0.50 & 1.00 
& 178.73 & 14.61 
& 178.90 & 14.79 
& 179.75 & 14.66 
& 180.85 & 14.55 \\
\midrule
0.75 & 0.50 
& 137.90 & 14.61 
& 139.16 & 14.82 
& 139.35 & 14.77 
& 140.52 & 14.74 \\

0.75 & 0.75 
& 136.57 & 15.06 
& 136.88 & 14.93 
& 137.14 & 14.87 
& 138.38 & 14.83 \\

0.75 & 1.00 
& 135.65 & 15.06 
& 135.35 & 15.01 
& 135.50 & 14.94 
& 136.98 & 14.89 \\
\midrule
1.00 & 0.50 
& 96.05 & 15.54 
& 97.40 & 15.29 
& 97.33 & 15.29 
& 98.77 & 15.31 \\

1.00 & 0.75 
& 95.46 & 15.54 
& 95.71 & 15.39 
& 95.56 & 15.39 
& 97.13 & 15.40 \\

1.00 & 1.00 
& 95.29 & 15.54 
& 94.47 & 15.46 
& 94.29 & 15.45 
& 96.05 & 15.46 \\
\bottomrule
\end{tabular}
\caption{Objective values and average prices under different pricing and dispatching policies under the demand scenario in which most riders travel from the central region to the peripheral regions.}
\label{tab:multi_region_pricing_moving_out}
\end{table}

Tables~\ref{tab:multi_region_pricing} and \ref{tab:multi_region_pricing_moving_out} report the objective values and average prices under different policies for the two directional-demand scenarios. In both tables, the last column, Zigzag, Dyn.\ (Loc.), corresponds to the same zigzag dispatching policy described in Section~\ref{subsec:simulation}, but implemented with localized supply and demand counts. Across both scenarios and all $(c_d, c_r)$ settings, the localized policy achieves the highest objective value. The improvement arises from using localized and proximity-weighted estimates of supply and demand, which allows the platform to more accurately capture the local supply-demand balance around each rider arrival location and to differentiate pricing and dispatching decisions across space. Moreover, the improvement is not driven by higher fares. The average price in the last column is often slightly lower than those under the other baselines. This suggests that the performance improvement comes from more accurate spatial targeting of pricing and dispatching adjustments.

\section{Conclusion} 
\label{sec:conclusion}

In this paper, we analyze a Markov model that approximates the rider queueing and pickup processes in ride-hailing/robotaxi platforms with a fixed number of drivers. We incorporate state-dependent service times to capture key spatial matching characteristics. The platform optimizes net profit, penalized by rider waiting, through adaptive pricing and dispatching decisions that are state-dependent. We show that, under mild assumptions, the optimal dispatching policy admits a closed-form expression with a zigzag structure---interpretable as a monotone threshold-type policy---and is tractable for price optimization thanks to the resulting closed-form stationary distribution. Building on this insight, we propose an efficient and scalable dynamic programming heuristic to approximate the optimal zigzag policy in more general settings.

Synthetic experiments demonstrate that our algorithm achieves near-optimal performance and excellent scalability. We further validate our approach in simulated ride-hailing dispatching and pricing environments that capture features in real-world spatial matching. We also adapt our algorithms to multi-region systems with spatial imbalances in demand and supply. The results show that the zigzag policy, computed without access to the simulation environment, performs at least as well as, and often better than, a constant-radius dispatching policy that is fine-tuned within the simulation environment. These findings underscore the practical promise of the proposed policy.

\bibliographystyle{informs2014} %
\bibliography{sample} %

\begin{APPENDICES}

\section{Additional Results}
\subsection{An Alternative Definition of Action Space}\label{subsec:alt_def_action_space}

In this section, we describe how the action space can be redefined to enable the use of uniformization in solving the objective in equation~\eqref{eq:objective_alt}. To transform the continuous-time system into an equivalent discrete-time framework, it is necessary to bound the transition rate in each state. At state $(l,m)$, we redefine the dispatching action as $(d_{a}(l,m), d_{c}(l,m))$, where $d_{a}(l,m)$ and $d_{c}(l,m)$ denote, respectively, the number of drivers to dispatch when the next event is a new rider joining the queue, or a trip completion, at state $(l,m)$. Since the number of dispatches cannot exceed the number of idle drivers and the number of riders in the queue at any state, the number of dispatches at state $(l,m)$ follows $d_{a}(l,m) \in \{0, 1, \ldots, \min\{L - l,\, m + 1\}\}$ and $d_{c}(l,m) \in \{0, 1, \ldots, \min\{L - l + 1,\, m\}\}$. If the next event is an effective rider arrival, the state transitions to $(l+d_{a}(l,m),m+1-d_{a}(l,m))$. Conversely, if the next event is a trip completion, the state transitions to $(l-1+d_{c}(l,m),m-d_{c}(l,m))$. Under this new action space, the transition rate at each state is bounded by $M_0:=\max_{l,m}\{l\mu_{l,m}\}+\Lambda$. Consequently, the Bellman optimality equation for the average reward MDP can be written as
\begin{align}\label{eq:bellman_opt}
    \max_{\substack{\lambda(l,m)\in[0,\Lambda], \\ 0\le d_{a}(l,m)\le\min\{L - l,\, m + 1\},\\ 0\le d_{c}(l,m)\le \min\{L - l + 1,\, m\}, \\ d_{a}(l,m),d_{c}(l,m)\in\mathbb Z_{\ge0}}} \Big\{ &
   \underbrace{\frac{\lambda(l,m)}{M_0}\!\Big(p_0 + w_p^r t_0 + p_1\!\big(\lambda(l,m)\big)t_0 + h\!\big(l+d_{a}(l,m),\, m+1-d_{a}(l,m)\big)-h(l,m)\Big)}_{\text{expected additional return under a new rider arrival}} \notag \\
    &+
    \underbrace{\frac{l \,\mu_{l,m}}{M_0}\,
    \Big(h\big(l-1+d_{c}(l,m),\, m-d_{c}(l,m)\big)-h(l,m)\Big)}_{\text{expected additional return under a trip completion}} -
    \underbrace{\frac{c_d\,l + c_r\,m}{M_0}}_{\text{penalty}} \Big\}=g.
\end{align}

In equation \eqref{eq:bellman_opt}, $g$ refers to the long-term average reward rate (revenue rate minus penalty rate) under the optimal policy. The bias term $h(l,m)$ refers to the difference in the long-run expected total reward gained by starting the process at state $(l,m)$ compared to starting the process with its steady-state distribution under an optimal policy. %
An optimal policy can be obtained by recovering the optimal solution of \eqref{eq:bellman_opt}. One can solve \eqref{eq:bellman_opt} using relative value iteration (see Chapter 8 of \citealt{Puterman1994MDP}) under an average reward setting.

\subsection{Sub-optimality of the Zigzag Policy}
\label{subsec:add_ex}

\smallskip

\begin{example}[Sub-optimality of the zigzag policy class]\label{ex:sub_opt_zigzag}
    Table \ref{tab:opt_nonzigzag_ex} gives an example where the optimal policy is non-zigzag. In this example, $p_1(\cdot)$ is constant. We set $c_d=c_r=0.1$.
    Table \ref{subtab:opt_nonzigzag_example_mu} gives the service rate $\mu$ at each state, which does not satisfy Assumption \ref{as:2}. Table \ref{subtab:opt_nonzigzag} shows the corresponding optimal dispatching policy, which does not follow a zigzag pattern. We compare this optimal dispatching policy with the optimal zigzag policy in Table \ref{subtab:opt_zigzag} to show why it is strictly better. We consider a coupling between systems under the two policies (we use system $(b)$ and system $(c)$ to represent the systems under the optimal dispatching policy and the optimal zigzag policy, respectively).
    We couple all random events until both systems reach state $(2,2)$ and make different dispatching decisions. We keep both systems the same effective request rate during coupling.
    At that state, the optimal dispatching policy chooses not to dispatch while the optimal zigzag policy dispatches a driver. %
    As a consequence, system $(c)$ dispatches and moves to state $(3,1)$, whereas system $(b)$ stays in state $(2,2)$. %
    One of the following three events can occur next: (1) a new rider arrival to both systems; (2) a trip completion in both systems; (3) a trip completion in system $(b)$ only. Event (3) is possible to happen because the trip completion rate at state $(2,2)$ ($2\mu_{2,2}=4$) exceeds one at state $(3,1)$ ($3\mu_{3,1}=3$). If event (1) or (2) occurs, the two systems return to the same state, and the coupling is maintained. If event (3) occurs, system $(b)$ completes an extra trip (relative to system $(c)$). At that point, system $(b)$ is at $(3,0)$ and system $(c)$ is at $(3,1)$. Then system $(b)$ will keep having the same throughput rate but smaller penalty rate until there is an additional trip completion for system $(c)$ compared to system $(b)$. After that occurs, systems $(b)$ and $(c)$ lie on the same anti-diagonal (i.e., both systems have the same value of $l+m$). Then either systems $b$ and $c$ are at the same state or system $b$ is at $(2,2)$ and system $c$ is at $(3,1)$. Both cases return to the situation we discussed before. Hence, during the coupling period, the throughput rate in both systems keeps the same, but system $(b)$ has a lower average penalty rate, which shows the optimal zigzag policy is indeed sub-optimal.

\end{example}

\begin{table}[htbp]
\small
    \begin{subtable}{.33\linewidth}
      \centering
        
        \begin{tabular}{c|ccccc}
         \backslashbox[10mm]{$l$}{$m$} & 0 & 1 & 2 & 3 & $\cdots$\\
         \hline 0 & 0 & 0 & 0 & 0 & $\cdots$\\
         1 & 1 & 1 & 2 & 2 & $\cdots$\\
         2 & 1 & 1 & 2 & 2 & $\cdots$\\
         3 & 1 & 1 & 2 & 2 & $\cdots$ \\
        \end{tabular}
        \caption{ Service rate $\mu$ at each state}
        \label{subtab:opt_nonzigzag_example_mu}
    \end{subtable}%
    \begin{subtable}{.33\linewidth}
      \centering
        
         \begin{tabular}{c|ccccc}
         \backslashbox[10mm]{$l$}{$m$} & 0 & 1 & 2 & 3 & $\cdots$\\
         \hline 0 & 0 & 1 & 1 & 1 & $\cdots$\\
         1 & 0 & 1 & 1 & 1 & $\cdots$\\
         2 & 0 & 1 & {\textbf{0}}\tikzmark{a} & 1& $\cdots$ \\
         3 & 0\tikzmark{c} & 0 & 0\tikzmark{b} & 0 & $\cdots$\\
        \end{tabular}
    \begin{tikzpicture}[overlay, remember picture, transform canvas={yshift=.25\baselineskip}]
    \draw [thick,-stealth,color=gray] ({pic cs:a}) [bend left] to ({pic cs:b}) node[above right, yshift=.25\baselineskip] {\scriptsize $\Lambda$};
    \draw [thick,-stealth,shorten <= 4pt,color=gray] ({pic cs:a}) [] to ({pic cs:c}) node[above, 
    xshift=.55\baselineskip,
    yshift=.25\baselineskip] {\scriptsize $2\mu_{2,2}$};
    \end{tikzpicture}
    \caption{Optimal dispatching policy}
    \label{subtab:opt_nonzigzag}
    \end{subtable} 
    \begin{subtable}{.33\linewidth}
      \centering
        
         \begin{tabular}{c|ccccc}
         \backslashbox[10mm]{$l$}{$m$} & 0 & 1 & 2 & 3 & $\cdots$\\
         \hline 0 & 0 & 1 & 1 & 1 & $\cdots$\\
         1 & 0 & 1 & 1 & 1 & $\cdots$\\
         2 & 0 & 1 & 1 & 1 & $\cdots$\\
         3 & 0\tikzmark{c1}  & \textbf{0}\tikzmark{a0} & 0\tikzmark{b1}  & 0 & $\cdots$
        \end{tabular}
    \begin{tikzpicture}
    [overlay, remember picture, transform canvas={yshift=.25\baselineskip},color=gray]

    \draw [thick,-stealth,shorten >= 3pt,color=gray] ({pic cs:a0}) [left] to ({pic cs:b1}) node[above, xshift=-.85\baselineskip] {\scriptsize $\Lambda$};
    \draw [thick,-stealth, shorten <= 5pt] ({pic cs:a0}) [] to ({pic cs:c1}) node[above, 
    xshift=.55\baselineskip] {\scriptsize $3\mu_{3,1}$};
    \end{tikzpicture}
    \caption{Optimal zigzag policy}
    \label{subtab:opt_zigzag}
    \end{subtable}  

    \caption{Example of the sub-optimality for zigzag dispatching policies. \normalfont{For all $m > 3$, column $m$ is identical to column $3$.}}
    \label{tab:opt_nonzigzag_ex}
\end{table}

\subsection{Structural Properties of Pricing}\label{subsec:pricing}
We explore some structural properties of optimal pricing under a zigzag path with finite path length. Given some zigzag dispatching policy with path $P$, let $\lambda_i:=\lambda(l_i,m_i)$ be the effective request rate at $(l_i,m_i), i\in \{1,\cdots,I\}$, and $v_i:=l_i\mu_{l_i,m_i}$ be the trip completion rate at that state. Following the notation in \eqref{eq:bellman_opt}, the Bellman optimality equation for $\lambda_i$ can be formulated as
\begin{equation}\label{eq:bellman_zigzag_pricing}
    \max_{\lambda_i\in [0,\Lambda]} \left\{\frac{\lambda_i}{M_0}\left(p_0+w_p^rt_0+p_1(\lambda_i)t_0+h(i+1)-h(i) \right) \right\} +\frac{v _i}{M_0}(h(i-1)-h(i)) -\frac{c_dl_i+c_rm_i}{M_0}=g,
\end{equation}
for $i=1,2,\cdots,I-1$. For $i=I$, we have to block the system by setting $\lambda_I=0$. Thus, we have
\begin{equation}\label{eq:bellman_zigzag_pricing_last}
\frac{v_I}{M_0}(h(I-1)-h(I)) -\frac{c_dl_I+c_rm_I}{ M_0}=g.
\end{equation}
Given an optimal zigzag dispatching policy, the optimal effective request rate and thus the optimal pricing policy can be computed by solving equations \eqref{eq:bellman_zigzag_pricing} and \eqref{eq:bellman_zigzag_pricing_last}. In what follows, we show that the resulting optimal price is non-decreasing in path index $i$ under some additional assumptions. We present a lemma first. Under $c_d=c_r\ge0$, let $P^*= ((l^\ast_1, m^\ast_1), (l^\ast_2, m^\ast_2), \ldots, (l^\ast_I, m^\ast_I))$ be the optimal path yielded by Theorem \ref{theorem:opt_policy_struct} and $v_i^*=l_i^\ast \mu_{l_i^\ast, m_i^\ast}$ be the corresponding trip completion rate.

\begin{lemma}\label{lemma:vi_mono}
    Suppose that $\mu_{l,m}$ is non-decreasing in $m$, Assumption \ref{as:2} holds and $c_d=c_r\geq 0$. Then $v_i^*$ is non-decreasing in $i$.
\end{lemma}
\begin{proof}{Proof of Lemma \ref{lemma:vi_mono}.}
 By Lemma \ref{lemma:serv_rate_ineq} and the construction of $P^*$ in Theorem \ref{theorem:opt_policy_struct}, we know $v_i^*$ is greater than or equal to the trip completion rate of all other states along the same anti-diagonal, which shows that $l_{i+1}^*\mu_{l_{i+1}^*,m_{i+1}^*}\geq l_i^*\mu_{l_i^*,m_i^*+1}$ holds for all $i$. Since $\mu_{l,m}$ is non-decreasing in $m$, we have
$$v_i^*=l_i^*\mu_{l_i^*,m_i^*}\leq l_i^*\mu_{l_i^*,m_i^*+1}\leq l_{i+1}^*\mu_{l_{i+1}^*,m_{i+1}^*}=v_{i+1}^*,
$$
which shows that $v_i^*$ is non-decreasing in $i$.\hfill\halmos
\end{proof}

\smallskip

Lemma \ref{lemma:vi_mono} establishes that, under mild conditions, the service rate $v_i^\ast$ is non-decreasing in $i$. In the following, we show that for any zigzag path $P$, whenever the corresponding $v_i$ is non-decreasing in $i$, the optimal effective request rate $\lambda_i$ is non-increasing in $i$ if $c_d=c_r=0$. Consequently, the associated optimal price is non-decreasing in $i$.

\begin{proposition}
[Surge Pricing]\label{theorem:surge_pricing}
Suppose $\lambda p_1(\lambda)$ is concave in $\lambda$. For any finite zigzag path $P$ with length $I$, if $c_d=c_r=0$ and $v_i$ is non-decreasing in $i\in \{1,\cdots,I\}$, then $\lambda_i$ is non-increasing in $i$ as well.
\end{proposition}

\begin{proof}{Proof of Proposition \ref{theorem:surge_pricing}.}
Let $\Delta_i:=h(i+1)-h(i)$ be the difference between the bias at $i+1$ and $i$. Let $G(\Delta):=\max_{\lambda \in[0, \Lambda]}\{\lambda(p_0+w_p^rt_0+p_1(\lambda)+\Delta)\}$. Then it is clear that $G(\Delta)$ is non-decreasing in $\Delta$ since any maximizer $\lambda$ must be greater or equal to 0. We assume $h(1)=0$ without loss of generality. Then equations \eqref{eq:bellman_zigzag_pricing} and \eqref{eq:bellman_zigzag_pricing_last} can be simplified as
\begin{equation}\label{eq:Bellman_pricing_simplified}
    G(\Delta_i)-v_i\Delta_{i-1}=g M_0,%
\end{equation}
for $i=1,\cdots, I$, where $v_1=0$ and $\Delta_I=-\infty$ (due to $\lambda_I=0$). As mentioned previously, bias $h(i)$ refers to the difference between the cumulative reward for a system that starts in state $(l_i,m_i)$ and one in which the stationary distribution determines the initial state. In the following, we show that $h(i)$ is strictly decreasing in $i$.

We consider the cumulative reward for a Markov chain starting at state $(l_i,m_i)$ and one starting at state $(l_{i+1},m_{i+1})$. We keep the effective request rate in the first chain always the same as that in the second chain unless the second chain is in state $(l_I,m_I)$ while the first chain not. After the two chains first meet at the same state, we let them evolve identically, so their future rewards are the same from that time onward.  Let $X_t^{(1)}:=(l_t^{(1)},m_t^{(1)})$ be the state of the first chain at time $t$ and $X_t^{(2)}:=(l_t^{(2)},m_t^{(2)})$ be the state of the second chain at time $t$. Let $\tau:=\inf\{t\ge 0: X_t^{(1)}=X_t^{(2)}\}$ be the coupling time, and let $
T_I^{(2)}:=\inf\{t\ge 0: X_t^{(2)}=(l_I,m_I)\}$ be the hitting time of state $I$ for the second chain. Since the zigzag path is finite and $(l_{i+1},m_{i+1})$ is strictly closer to state $I$ than $(l_i,m_i)$, there is a positive-probability event on which the second chain reaches state $(l_I,m_I)$ before the two chains meet: $\mathbb{P}(T_I^{(2)}<\tau)>0$. On this event, the second chain must set $\lambda_I=0$ and therefore earns no further revenue after time $T_I^{(2)}$ before next event happens. By contrast, the first chain is still in some state $j<I$, where a strictly positive price can be chosen, yielding a strictly positive revenue rate. Since $c_d=c_r=0$, this creates a strictly positive reward gap on the event $\{T_I^{(2)}<\tau\}$. This implies $\Delta_i=h(i+1)-h(i)<0$ for $i=1,\cdots,I$.

Next, we show that $\Delta_i$ is strictly decreasing in $i$. We first notice $\Delta_1< \Delta_0:=0$. By induction, for any $k=1,\cdots,I-1$, given $\Delta_k< \Delta_{k-1}< \cdots< \Delta_1 < 0$, we need to show $\Delta_{k+1}< \Delta_k$. By setting $i=k$ and $i=k+1$ in equation \eqref{eq:Bellman_pricing_simplified} and taking difference on both sides, we get the following equation after cancelling out common term $gM_0$:
\begin{equation}\label{eq:bellman_induction_eq}
    G(\Delta_k)-G(\Delta_{k+1})-v_k\Delta_{k-1}+v_{k+1}\Delta_{k}=0.%
\end{equation}
Since $\Delta_k< \Delta_{k-1}< 0$ and $v_{k+1},v_k>0$, equation \eqref{eq:bellman_induction_eq} yields the following inequality:
\begin{equation}\label{eq:middle_terms_induction}
    G(\Delta_k)-G(\Delta_{k+1})=v_k\Delta_{k-1}-v_{k+1}\Delta_{k}\geq v_{k}(\Delta_{k-1}-\Delta_{k})> 0
\end{equation}
Since $G(\Delta)$ is non-decreasing in $\Delta$, $G(\Delta_k)> G(\Delta_{k+1})$ leads to $\Delta_{k+1}< \Delta_k$, which completes the induction proof. 

Finally, we show $\lambda_i$ is non-increasing in $i$. Let $F(\lambda, \Delta):=\lambda\left(p_0+w_p^r t_0+p_1(\lambda)+\Delta\right)$. For each $i$, the Bellman equation chooses $\lambda_i \in \arg \max _{\lambda \in[0, \Lambda]} F\left(\lambda, \Delta_i\right)$. Since we have proved that $\Delta_i$ is strictly decreasing in $i$, fix any $i \in\{1, \ldots, I-1\}$. Suppose for contradiction that $\lambda_{i+1}>\lambda_i$. By optimality of $\lambda_i$ for $\Delta_i$ and of $\lambda_{i+1}$ for $\Delta_{i+1}$, we have
$
F\left(\lambda_i, \Delta_i\right) \geq F\left(\lambda_{i+1}, \Delta_i\right)$ and $F\left(\lambda_{i+1}, \Delta_{i+1}\right) \geq F\left(\lambda_i, \Delta_{i+1}\right)
$. Adding these two inequalities gives
$$
\left(\lambda_i-\lambda_{i+1}\right)\left(\Delta_i-\Delta_{i+1}\right) \geq 0 .
$$
But $\lambda_i-\lambda_{i+1}<0$ and $\Delta_i-\Delta_{i+1}>0$, so the left-hand side is strictly negative, which leads to a contradiction. Hence we know $\lambda_{i+1} \leq \lambda_i$ holds for all $i$. Thus, $\lambda_i$ is non-increasing in $i$. \hfill\halmos
\end{proof}

\medskip

Proposition~\ref{theorem:surge_pricing}, together with Lemma~\ref{lemma:vi_mono}, implies that when the rider queue has finite capacity, the penalty coefficients satisfy $c_d=c_r=0$, and $\mu_{l,m}$ is nondecreasing in $m$, the optimal pricing policy characterized in Theorem~1 is nondecreasing along the zigzag path. The intuition is as follows. When the path index $i$ is small, the trip completion rate is relatively low, so a lower price induces a higher effective request rate, which not only increases the current revenue rate but also helps move the system toward states with larger $i$, where trips are completed more quickly. By contrast, when $i$ is already large, the system is able to complete trips rapidly, and additional arrivals mainly hasten the transition toward the terminal state, where the platform earns no revenue. The platform therefore sets a higher price in such states to discourage arrivals and reduce the likelihood of reaching that terminal state. This monotone pricing structure resembles the \emph{surge pricing} commonly observed on ride-hailing platforms: when the number of nearby idle drivers is low, the platform increases rider prices to better balance supply and demand. Proposition \ref{theorem:surge_pricing} extends this intuition to a two-sided spatial queueing system in which riders may also wait in queue.

\medskip

The following proposition shows that when $c_d=c_r\geq 0$, monotonicity of the optimal pricing policy can still be established if $v_i$ is strictly decreasing in $i$. In this case, as $i$ increases, the penalty rate rises while the trip completion rate falls. The platform therefore has an incentive to keep the system in states with smaller $i$ in order to achieve a higher long-run objective, which implies that the optimal price should increase along the path. This situation may arise in a spatial setting when the zigzag path moves continuously downward and $\mu_{l,m}$ decreases sharply with $l$.

\begin{proposition}
\label{corollary:vi=iv_1}
    Suppose $\lambda p_1(\lambda)$ is concave in $\lambda$ and $c_d=c_r \geq 0$. For any finite zigzag path $P$ with length $I$, if $v_i$ is strictly decreasing in $i\in \{2,\cdots,I\}$, then $\lambda_i$ is non-increasing in $i$.
\end{proposition}
\begin{proof}{Proof of Proposition \ref{corollary:vi=iv_1}.}
Similar to \eqref{eq:Bellman_pricing_simplified}, equation~\eqref{eq:bellman_zigzag_pricing} can be rewritten as
\begin{equation}\label{eq:Bellman_simplified_penalty}
    G(\Delta_i)-v_i\Delta_{i-1}=gM_0+(c_d l_i+c_r m_i),
\end{equation}
for $i=1,\ldots,I$. We again show that $\Delta_i$ is strictly decreasing in $i$. First, we set $\Delta_{I-1}>\Delta_I:=-\infty$. By a coupling argument, one can similarly show that $\Delta_i=h(i+1)-h(i)<0$. Now suppose, for some $k\in\{2,\ldots,I-1\}$, that
\[
\Delta_k>\Delta_{k+1}>\cdots>\Delta_I.
\]
To complete the induction, it suffices to prove that $\Delta_{k-1}>\Delta_k$. Setting $i=k$ and $i=k+1$ in \eqref{eq:Bellman_simplified_penalty} and taking the difference, after canceling the common term $gM_0$, yields
\begin{equation}\label{eq:bellman_induction_eq_penalty}
    G(\Delta_k)-G(\Delta_{k+1})-v_k\Delta_{k-1}+v_{k+1}\Delta_k
    =(c_d l_k+c_r m_k)-(c_d l_{k+1}+c_r m_{k+1})\leq 0,
\end{equation}
where the inequality follows from $l_k+m_k\leq l_{k+1}+m_{k+1}$ and $c_d=c_r$. Since $\Delta_k<\Delta_{k-1}<0$, equation~\eqref{eq:bellman_induction_eq_penalty} implies
\begin{equation}\label{eq:middle_terms_induction_penalty}
    v_k(\Delta_{k-1}-\Delta_k)
    > v_k\Delta_{k-1}-v_{k+1}\Delta_k
    \geq G(\Delta_k)-G(\Delta_{k+1})
    \geq 0,
\end{equation}
which in turn yields $\Delta_{k-1}>\Delta_k$. Therefore, if $v_i$ is nonincreasing in $i$, then $\Delta_i$ is strictly decreasing in $i$. By the same argument as in the proof of Proposition~\ref{theorem:surge_pricing}, it follows that $\lambda_i$ is non-increasing in $i$.  \hfill\halmos
\end{proof}

\subsection{Additional Performance Metrics in Synthetic Experiments}\label{sec:add_perf_metric_sync}

We provide some additional performance metrics including average price for each rider, revenue rate, average pickup time, and average queueing time under the experiment setting in Section \ref{subsec:synthetic_exp}. Tables \ref{tab:zigzag_metrics_20} and \ref{tab:zigzag_metrics} show the relevant metrics under scenarios $(L,\Lambda,M)=(20,8,10)$ and $(L,\Lambda,M)=(100,40,50)$, respectively. When $(L,\Lambda,M)=(20,8,10)$, all metrics for the dynamic-pricing zigzag policy are close to those obtained by value iteration. By contrast, the greedy policy yields significantly larger pickup times and shorter queueing times compared to other policies. Under the zigzag policy, static pricing results in a higher average price and comparable pickup time but leads to longer queueing times compared to dynamic pricing. This is because static pricing lacks the ability to adaptively dampen demand as congestion builds, causing the system to rely more heavily on queueing.

\begin{table}[htbp]
\SingleSpacedXI
\centering
\small
\begin{tabular}{cc cccc cccc}
\toprule
           &            & \multicolumn{4}{c}{\textbf{Value Iteration}} &
             \multicolumn{4}{c}{\textbf{Greedy, Dyn.}} \\
\cmidrule(lr){3-6}\cmidrule(lr){7-10}
$c_d$ & $c_r$ & Revenue & Price & Pickup & Queue
      & Revenue & Price & Pickup & Queue \\
\midrule
0.50 & 0.50 & 28.30 & 10.86 & 1.36 & 1.53  & 27.41 & 10.93 & 2.55 & 1.65 \\
0.50 & 0.75 & 27.59 & 11.00 & 1.47 & 1.18  & 25.70 & 11.36 & 2.79 & 0.89 \\
0.50 & 1.00 & 26.86 & 11.12 & 1.67 & 0.86  & 23.73 & 11.60 & 2.77 & 0.26 \\
\midrule
0.75 & 0.50 & 27.74 & 11.01 & 1.23 & 1.50  & 24.84 & 11.43 & 2.55 & 0.83 \\
0.75 & 0.75 & 26.90 & 11.16 & 1.40 & 1.07  & 22.31 & 11.75 & 2.38 & 0.07 \\
0.75 & 1.00 & 26.14 & 11.27 & 1.49 & 0.84  & 22.04 & 11.75 & 2.33 & 0.02 \\
\midrule
1.00 & 0.50 & 26.72 & 11.23 & 1.13 & 1.42  & 21.00 & 11.95 & 2.09 & 0.00 \\
1.00 & 0.75 & 25.75 & 11.38 & 1.28 & 1.00  & 20.98 & 11.95 & 2.08 & 0.00 \\
1.00 & 1.00 & 24.93 & 11.50 & 1.42 & 0.73  & 20.98 & 11.95 & 2.08 & 0.00 \\
\end{tabular}

\begin{tabular}{cc cccc cccc}
\toprule
           &            & \multicolumn{4}{c}{\textbf{Zigzag, Dyn.}} &
             \multicolumn{4}{c}{\textbf{Zigzag, Sta.}} \\
\cmidrule(lr){3-6}\cmidrule(lr){7-10}
$c_d$ & $c_r$ & Revenue & Price & Pickup & Queue
      & Revenue & Price & Pickup & Queue \\
\midrule
0.50 & 0.50 & 28.30 & 10.86 & 1.36 & 1.53  & 27.75 & 10.97 & 1.38 & 1.59 \\
0.50 & 0.75 & 27.60 & 11.00 & 1.44 & 1.20  & 27.16 & 11.14 & 1.45 & 1.32 \\
0.50 & 1.00 & 26.87 & 11.12 & 1.65 & 0.88  & 26.46 & 11.29 & 1.62 & 1.02 \\
\midrule
0.75 & 0.50 & 27.75 & 11.01 & 1.26 & 1.45  & 27.17 & 11.13 & 1.27 & 1.50 \\
0.75 & 0.75 & 26.90 & 11.16 & 1.40 & 1.07  & 26.40 & 11.31 & 1.40 & 1.16 \\
0.75 & 1.00 & 26.14 & 11.27 & 1.49 & 0.84  & 25.62 & 11.45 & 1.48 & 0.93 \\
\midrule
1.00 & 0.50 & 26.73 & 11.23 & 1.15 & 1.39  & 26.16 & 11.35 & 1.17 & 1.42 \\
1.00 & 0.75 & 25.75 & 11.38 & 1.29 & 0.99  & 25.23 & 11.52 & 1.29 & 1.06 \\
1.00 & 1.00 & 24.93 & 11.50 & 1.42 & 0.73  & 24.37 & 11.65 & 1.41 & 0.80 \\
\bottomrule
\end{tabular}
\caption{Additional performance metrics of different algorithms under parameters $(L,\Lambda,M)=(20,8,10)$.}
\label{tab:zigzag_metrics_20}
\end{table}

\begin{table}[htbp]
\SingleSpacedXI
\centering
\small
\begin{tabular}{cc cccc cccc}
\toprule
           &            & \multicolumn{4}{c}{\textbf{Value Iteration}} &
             \multicolumn{4}{c}{\textbf{Greedy, Dyn.}} \\
\cmidrule(lr){3-6}\cmidrule(lr){7-10}
$c_d$ & $c_r$ & Revenue & Price & Pickup & Queue
      & Revenue & Price & Pickup & Queue \\
\midrule
0.50 & 0.50 & 164.24 &  9.56 & 0.41 & 2.62  & 160.99 & 10.25 & 1.15 & 1.45 \\
0.50 & 0.75 & 164.59 &  9.59 & 0.47 & 2.56  & 150.34 & 10.59 & 1.34 & 0.73 \\
0.50 & 1.00 & 165.01 &  9.61 & 0.47 & 2.54  & 139.13 & 11.00 & 1.33 & 0.00 \\
\midrule
0.75 & 0.50 & 164.34 &  9.58 & 0.38 & 2.66  & 160.95 & 10.25 & 1.15 & 1.45 \\
0.75 & 0.75 & 164.79 &  9.59 & 0.41 & 2.60  & 137.48 & 11.09 & 1.22 & 0.00 \\
0.75 & 1.00 & 165.08 &  9.61 & 0.46 & 2.54  & 137.48 & 11.09 & 1.22 & 0.00 \\
\midrule
1.00 & 0.50 & 163.89 &  9.65 & 0.34 & 2.79  & 134.07 & 11.24 & 1.09 & 0.00 \\
1.00 & 0.75 & 164.74 &  9.60 & 0.38 & 2.65  & 134.07 & 11.24 & 1.09 & 0.00 \\
1.00 & 1.00 & 165.23 &  9.61 & 0.41 & 2.58  & 134.07 & 11.24 & 1.09 & 0.00 \\
\end{tabular}

\begin{tabular}{cc cccc cccc}
\toprule
           &            & \multicolumn{4}{c}{\textbf{Zigzag, Dyn.}} &
             \multicolumn{4}{c}{\textbf{Zigzag, Sta.}} \\
\cmidrule(lr){3-6}\cmidrule(lr){7-10}
$c_d$ & $c_r$ & Revenue & Price & Pickup & Queue
      & Revenue & Price & Pickup & Queue \\
\midrule
0.50 & 0.50 & 163.11 & 10.17 & 0.54 & 0.71  & 161.82 & 10.27 & 0.52 & 0.78 \\
0.50 & 0.75 & 160.93 & 10.26 & 0.58 & 0.52  & 159.38 & 10.38 & 0.56 & 0.58 \\
0.50 & 1.00 & 159.27 & 10.32 & 0.62 & 0.42  & 157.49 & 10.46 & 0.59 & 0.47 \\
\midrule
0.75 & 0.50 & 161.13 & 10.27 & 0.51 & 0.62  & 159.70 & 10.37 & 0.50 & 0.67 \\
0.75 & 0.75 & 158.75 & 10.36 & 0.55 & 0.46  & 157.07 & 10.47 & 0.53 & 0.50 \\
0.75 & 1.00 & 157.03 & 10.43 & 0.58 & 0.38  & 155.13 & 10.55 & 0.56 & 0.41 \\
\midrule
1.00 & 0.50 & 157.30 & 10.44 & 0.47 & 0.52  & 156.03 & 10.51 & 0.47 & 0.55 \\
1.00 & 0.75 & 154.81 & 10.54 & 0.51 & 0.39  & 153.36 & 10.62 & 0.50 & 0.42 \\
1.00 & 1.00 & 152.99 & 10.60 & 0.54 & 0.33  & 151.38 & 10.69 & 0.53 & 0.35 \\
\bottomrule
\end{tabular}
\caption{Additional performance metrics of different algorithms under parameters $(L,\Lambda,M)=(100,40,50)$.}
\label{tab:zigzag_metrics}
\end{table}

\subsection{Performance Evaluation of Algorithm 1}\label{sec:perf_evaluation_alg_1}
Algorithm~1 is based on a dominance heuristic along zigzag paths: when two candidate paths terminate at the same state, the one with a higher static objective value is used as the preferred prefix when extending the path. This section evaluates how reliably this heuristic leads to optimal or near-optimal performance in instances where the optimal policy can be computed. We first formalize and prove the exact dominance condition discussed in Section \ref{sec:zigzag_policy}:

\begin{proposition}\label{prop:path_comparison}
Consider two paths $P=\{(l_{i},m_{i}),i=1,\cdots,I\}$ and $P^\prime=\{(l^\prime_{i},m^\prime_{i}),i=1,\cdots,I\}$. If $v_i\leq v_i^\prime$ and $c_dl_i+c_rm_i\geq c_dl_i^\prime+c_rm_i^\prime$ holds for all $i$, then $P^\prime$ dominates $P$. In other words, for any future extensions of same states on both $P$ and $P^\prime$, the objective value associated with $P^\prime$ under the corresponding optimal pricing policy is always greater than or equal to one associated with $P$. 
\end{proposition}

\begin{proof}{Proof of Proposition \ref{prop:path_comparison}.}

Fix arbitrary future extensions of $P$ and $P'$. Let $\lambda^*$ be an optimal pricing policy for $P$ together with its extension. We construct a feasible pricing policy $\lambda$ for $P'$ and its extension such that the effective request rate in system 2 (the system evolving along $P'$ and its extension) always matches that in system 1 (the system evolving along $P$ and its extension), even when the two systems are in different states. This construction is feasible because system 2 can never be at a state with a larger path index than system 1 at any time, as we show below. Suppose both systems start from the same state, and couple them using the same realization of all exogenous randomness, including arrivals and trip completions. Under this coupling, since the effective request rates are identical at all times, the two systems also earn the same revenue rate at all times.

We next compare the state evolution of the two systems along their respective paths. Whenever both systems are at index $i$, system 1 is at state $(l_i,m_i)$ and system 2 is at state $(l_i',m_i')$. Since $v_i' \geq v_i$, system 2 has a trip completion rate weakly larger than that of system 1 at the same index. Under the coupling, a trip completion cannot occur in system 1 without also occurring in system 2. Therefore, starting from any time at which the two systems share the same index $i$, the next state index reached by system 2 can never exceed that reached by system 1. Repeating this argument each time the two systems meet at the same index, we conclude by induction that system 2 can never have a larger path index than system 1 along the coupled sample path.

Now compare the penalty accumulated at each index $i$. By assumption, $c_d l_i + c_r m_i \geq c_d l_i' + c_r m_i'$, and the penalty rate under both $P$ and $P'$ is increasing in the path index $i$. Since system 2 always has an index weakly smaller than that of system 1, its penalty rate is always weakly smaller than that of system 1. Because the revenue rates are identical in the two systems, the objective value in system 2 is weakly greater than that in system 1. Therefore, under the corresponding optimal pricing policies, the objective value associated with $P'$ is always weakly greater than that associated with $P$. \hfill$\halmos$

\end{proof}

\medskip

The dominance condition in Proposition~\ref{prop:path_comparison} is, however, quite restrictive. Along two distinct zigzag paths, the inequalities $v_i' \geq v_i$ and $c_d l_i + c_r m_i \geq c_d l_i' + c_r m_i'$ must hold at every index $i$, a requirement that is rarely satisfied in practice. As a result, this exact dominance rule alone is not sufficient to limit the number of candidate prefixes retained at each terminal state, and Algorithm~\ref{alg:DP_zigzag} must still rely on a heuristic comparison of static objective values when deciding which prefix to extend. This motivates a more direct numerical assessment of how reliably the heuristic delivers optimal or near-optimal performance in settings where the true optimal policy can still be computed, which we present next.

We consider number of vehicles $L=5,6,\ldots,30$. For each $L$, we generate 20 independent instances. In each instance, we use the same way described in Section \ref{subsec:synthetic_exp} to estimate service rate $\mu$ at each state. We sample the arrival rate as $\Lambda \sim \operatorname{Unif}(0.2L,0.6L)$ and the penalty coefficients $c_d,c_r \sim \operatorname{Unif}(0.5,1)$ independently. For every instance, we compute the objective achieved by Algorithm \ref{alg:DP_zigzag} and the optimal objective value through value iteration (which is tractable for this range of $L$).

\begin{figure}[htbp]
\centering
\begin{tikzpicture}
\begin{groupplot}[
    group style={
        group size=3 by 1,
        horizontal sep=1.8cm
    },
    width=0.32\textwidth,
    height=0.32\textwidth,
    axis lines=left,
    xmin=4,
    xmax=31,
    ymin=0,
    xtick distance=5,
    tick label style={font=\small},
    label style={font=\small},
]

\nextgroupplot[
title={(a)},
    xlabel={$L$},
    ylabel={Optimality Frequency (\%)},
    xmin=4,
    xmax=31,
    ymax=100,
]
\addplot+[
    ycomb,
    blue,
    very thin,
    mark=*,
    mark size=0.8pt,
] table[
    col sep=comma,
    x=x,
    y=y
]{figures/opt_pct_match.csv};

\nextgroupplot[
title={(b)},
    xlabel={$L$},
    ylabel={Average Performance Ratio},
    ymin=0.998,
    ymax=1,
    ytick={0.998,0.999,1},
    yticklabels={0.998,0.999,1},
    legend style={at={(0.97,0.03)}, anchor=south east}
]
\addplot+[
    ycomb,
    blue,
    very thin,
    mark=*,
    mark size=0.8pt,
] table[
    col sep=comma,
    x=x,
    y=y
]{figures/avg_ratio_match.csv};

\nextgroupplot[
title={(c)},
    xlabel={$L$},
    ylabel={Worst Performance Ratio},
    xmin=4,
    xmax=31,
    ymin=0.998,
    ymax=1,
    ytick={0.998,0.999,1},
    yticklabels={0.998,0.999,1},
]
\addplot+[
    ycomb,
    blue,
    very thin,
    mark=*,
    mark size=0.8pt,
] table[
    col sep=comma,
    x=x,
    y=y
]{figures/worst_ratio_match.csv};

\end{groupplot}
\end{tikzpicture}
\caption{Performance metrics as functions of $L$. From left to right, the panels report the optimal percentage, the average ratio, and the worst-case ratio.}
\label{fig:synthetic_results_three_panel}
\end{figure}

Figure~\ref{fig:synthetic_results_three_panel}(a) plots the optimality frequency, defined as the fraction of instances in which Algorithm~\ref{alg:DP_zigzag} attains the exact optimal value based on the full dispatching action space defined in Appendix~\ref{subsec:alt_def_action_space}, as a function of $L$. We observe a relatively high optimality frequency (above $60\%$) for small instances with $L<10$. As $L$ increases, however, the number of comparisons among candidate paths grows exponentially, which increases the likelihood that the heuristic dominance check fails. As a result, preserving exact optimality throughout the entire path construction process becomes increasingly difficult when $L$ is large. Despite this decline in exact-optimality frequency, Figures~\ref{fig:synthetic_results_three_panel}(b) and~\ref{fig:synthetic_results_three_panel}(c) show that the objective values produced by Algorithm~\ref{alg:DP_zigzag} remain extremely close to optimal. The average performance ratio, defined as the value produced by Algorithm~\ref{alg:DP_zigzag} divided by the exact optimal value, stays within a very narrow band around $1$ for all $L$. Moreover, for each $L$, the worst-case ratio across the 20 instances also remains uniformly high. In particular, the worst observed ratio over all tested instances exceeds $0.998$, indicating that Algorithm~\ref{alg:DP_zigzag} achieves at least $99.8\%$ of the optimal objective in every instance we tested.

\end{APPENDICES}

\ACKNOWLEDGMENT{The authors thank the support of an Amazon Science Faculty Research Award
through the UW-Amazon Science Hub.}

\ECSwitch

\ECHead{\centering{Online Appendix}}

\section{Proofs}

\smallskip

\begin{proof}{Proof of Lemma \ref{lemma:equivalent_penalty}.}
First, we notice that the sum of the expected number of drivers in service and idle drivers equals the total number of drivers in the system, so we have $L_o^d(\lambda,\phi)=L-L_s^d(\lambda,\phi)$. Let $L_t^d(\lambda,\phi)$ be the expected number of riders on a trip under policy $(\lambda,\phi)$. Since the total number of riders in service includes both those waiting to be picked up and those currently on a trip, we have $L_s^d(\lambda,\phi)=L_p^r(\lambda,\phi)+L_t^d(\lambda,\phi)$. We consider an observation period $T$. Over this period, $L_t^d(\lambda, \phi)T$ is the total on-trip time for all riders. $L_t^d(\lambda,\phi)T/t_0$ is the expected number of riders served by the platform over period $T$. As $T \rightarrow \infty$, by the law of large numbers, the average on-trip time converges to $t_0$. Thus, the average number of riders served by the platform per time unit can be computed by $\lim_{T\rightarrow\infty} L_t^d(\lambda,\phi)T/(t_0T)=L_t^d(\lambda,\phi)/t_0$. On the other hand, if we increase the base fare $p_0$ by 1 unit, the objective $\tilde{\mathcal{R}}(p_0+1,\lambda, \phi)$ will increase by the average number of riders served per time unit according to equation \eqref{eq:objective}. %
Thus, if the platform increases the base fare by $w_p^rt_0$, we have $\mathcal{R}(p_0+w_p^rt_0,\lambda,\phi)=\mathcal{R}(p_0,\lambda,\phi)+w_p^rt_0L_t^d(\lambda,\phi)/t_0$. Consequently, equation \eqref{eq:objective} can be rewritten as
\begin{align*}
    &\mathcal{R}(p_0,\lambda,\phi)-w_s^dL_s^d(\lambda,\phi)-w_q^rL_q^r(\lambda,\phi)-w_p^rL_p^r(\lambda,\phi)-w_o^dL_o^d(\lambda,\phi) \\
    =&\mathcal{R}(p_0,\lambda,\phi)-w_s^dL_s^d(\lambda,\phi)-w_q^rL_q^r(\lambda,\phi)-w_p^r\left(L_s^d(\lambda,\phi)-L_t^d(\lambda,\phi)\right)-w_o^d(L-L_s^d(\lambda,\phi))\\
    =&\mathcal{R}(p_0,\lambda,\phi)+\frac{w_p^rt_0L_t^d(\lambda,\phi)}{t_0}-(w_s^d+w_p^r-w_o^d)L_s^d(\lambda,\phi)-w_q^rL_q^r(\lambda,\phi)-w_o^dL \\
    =&\mathcal{R}(p_0+w_p^rt_0,\lambda,\phi)-(w_s^d+w_p^r-w_o^d)L_s^d(\lambda,\phi)-w_q^rL_q^r(\lambda,\phi)-w_o^dL,%
\end{align*}
which proves the result.%
\hfill \halmos
\end{proof}

\medskip

\begin{proof}{Proof of Lemma \ref{lemma:serv_rate_ineq}.}
We start with the proof of condition (1). Let $(l^\prime,m^\prime)$ be a type 1 state. To show that $(l,m)$ is also a type 1 state for any $l^\prime \leq l \leq L-1$ and $1\leq m\leq m^\prime$, %
it is sufficient to show that $(l^\prime+1,m^\prime)$ and $(l^\prime,m^\prime-1)$ are both type 1 state. 
Then condition (1) can be proved by recursively applying this result. We first show that $(l^\prime+1,m^\prime)$ is a type 1 state. According to Assumption \ref{as:2}, we have
\begin{equation}\label{eq:submod}
    (l^\prime+1)\mu_{l^\prime+1,m^\prime}-(l^\prime+1)\mu_{l^\prime+1,m^\prime-1}\geq l^\prime\mu_{l^\prime,m^\prime}-l^\prime\mu_{l^\prime,m^\prime-1},
\end{equation}
\begin{equation}\label{eq:submod2}
\begin{aligned}
    (l^\prime+1)\mu_{l^\prime+1,m^\prime-1}-(l^\prime+2)\mu_{l^\prime+2,m^\prime-1}&\geq l^\prime\mu_{l^\prime,m^\prime-1}-(l^\prime+1)\mu_{l^\prime+1,m^\prime-1}.
    \end{aligned}
\end{equation}
Inequality \eqref{eq:submod} is from condition (1) in Assumption \ref{as:2}, and inequality \eqref{eq:submod2} is from condition (2) in Assumption \ref{as:2}. Summing equation \eqref{eq:submod} with equation \eqref{eq:submod2} on both sides yields
\begin{equation}
    \label{eq:sum_sub}
    (l^\prime+1)\mu_{l^\prime+1,m^\prime}-(l^\prime+2)\mu_{l^\prime+2,m^\prime-1}\geq l^\prime\mu_{l^\prime,m^\prime}-(l^\prime+1)\mu_{l^\prime+1,m^\prime-1}>0,
\end{equation}
where the second inequality is due to the fact that $(l^\prime,m^\prime)$ is a type 1 state which satisfies $l^\prime \mu_{l^\prime, m^\prime} >(l^\prime+1) \mu_{l^\prime+1, m^\prime-1}$. Equation \eqref{eq:sum_sub} yields $(l^\prime+1)\mu_{l^\prime+1,m^\prime}>(l^\prime+2)\mu_{l^\prime+2,m^\prime-1}$, which shows that $(l^\prime+1,m^\prime)$ is a type 1 state by definition. To show $(l^\prime,m^\prime-1)$ is also a type 1 state, similarly, Assumption \ref{as:2} gives
\begin{equation}\label{eq:submod3}
\begin{aligned}
    l^\prime\mu_{l^\prime,m^\prime-1}-l^\prime\mu_{l^\prime,m^\prime-2}&\geq l^\prime\mu_{l^\prime,m^\prime}-l^\prime\mu_{l^\prime,m^\prime-1},
    \end{aligned}
\end{equation}
\begin{equation}\label{eq:submod4}
    l^\prime\mu_{l^\prime,m^\prime-2}-(l^\prime+1)\mu_{l^\prime+1,m^\prime-2}\geq l^\prime\mu_{l^\prime,m^\prime-1}-(l^\prime+1)\mu_{l^\prime+1,m^\prime-1}.
\end{equation}
Summing equation \eqref{eq:submod3} with equation \eqref{eq:submod4} on both sides, we get
\begin{equation}\notag%
    l^\prime\mu_{l^\prime,m^\prime-1}-(l^\prime+1)\mu_{l^\prime+1,m^\prime-2}\geq l^\prime\mu_{l^\prime,m^\prime}-(l^\prime+1)\mu_{l^\prime+1,m^\prime-1}>0,
\end{equation}
which shows that $(l^\prime,m^\prime-1)$ is a type 1 state.

So far, we have proved condition (1). Now we consider the case when $(l^\prime,m^\prime)$ is a type 2 state. If there exists some type 1 state $(l_1,m_1)$ such that $l_1\leq l^\prime$ and $m_1\geq m^\prime$, according to condition 1, we know $(l^\prime,m^\prime)$ must be a type 1 state as well, which leads to contradiction. Hence, we know condition (2) holds. \hfill\halmos
\end{proof}

\medskip

\begin{proof}{Proof of Theorem \ref{theorem:opt_policy_struct}.}
\renewcommand\labelenumi{(\theenumi)}
\setlist[enumerate]{itemsep=-0.4mm}
\setlist[itemize]{itemsep=-0.4mm}

We prove the theorem under three cases $c_d=c_r$, $c_d>c_r$, and $c_d<c_r$. In all cases, for an arbitrary dispatching policy, we construct a new policy that satisfies the optimality properties stated in the theorem. We then show that this constructed policy achieves a weakly higher objective value than the original policy.

\smallskip

\noindent\underline{\textbf{Case 1: $c_d = c_r$}}. Let $\phi^*$ be the dispatching policy such that $\phi^*(l,m)=0$ for all type 1 states and $\phi^*(l,m)=1$ for all type 2 states and $\lambda^\ast$ the optimal dynamic pricing policy under $\phi^*$. We want to show $\Tilde{\mathcal{R}}(\phi^*,\lambda^*)\geq \Tilde{\mathcal{R}}(\phi,\lambda)$ for an arbitrary policy $(\phi,\lambda)$. We consider a coupling of two systems, where system 1 uses policy $(\phi^*,\lambda^*)$ and system 2 uses the policy $(\phi,\lambda)$. This coupling aligns all random events in both systems, including arrivals, trip completions, and dispatches. Instead of using the optimal pricing policy $\lambda^*$ for system 1, we apply a dependent pricing policy $\lambda^{\mathrm{D}}$ that is dependent on system 2. Specifically, $\lambda^{\mathrm{D}}$ is constructed to keep the effective request rate in system 1 always the same as that in system 2, even if the two systems are in different states. In other words, for any state transition in system 2, the effective request rate in system 1 is immediately adjusted to the same as the new effective request rate in system 2. %
We know $\Tilde{\mathcal{R}}(\phi^*,\lambda^*)\geq \Tilde{\mathcal{R}}(\phi^*,\lambda^{\mathrm{D}})$. Thus, it is sufficient to show $\Tilde{\mathcal{R}}(\phi^*,\lambda^{\mathrm{D}})\geq \Tilde{\mathcal{R}}(\phi,\lambda)$. We demonstrate this in the following.

Let \(\{X_t\}_{t=1}^\infty\) and \(\{Y_t\}_{t=1}^\infty\) denote the coupled Markov chain, where each step $t$ corresponds to a transition in the joint process. At step \(t\), the state of system 1 is given by \(X_t = (l_t, m_t)\), and the state of system 2 is given by \(Y_t = (l_t', m_t')\). Without loss of generality, we assume that both systems start with the same number of drivers in service and riders in the queue, that is, $X_1=Y_1$. If a random event occurs in either system at step $t$, the states of systems 1 and 2 transition to $X_{t+1}$ and $Y_{t+1}$ after possible dispatchings, respectively, based on their dispatching policies. After that, the coupled chain proceeds to step $t+1$.

Note that any arrival events will occur in both chains at the same time in the coupled system because the effective request rate in system 1 is always set to be the same as that in system 2. However, a trip completion may occur in only one of the chains when the service completion rate in one system is larger than the other.

In the following proof, we show that under coupling, there are three possible ``phases’’, which are mutually exclusive, for the pair $(X_t, Y_t)$ at step $t$. The coupled process loops over these phases. Moreover, we will show that, at each phase, the penalty rate in system 1 is consistently less than or equal to one in system 2. 

\smallskip

\textbf{Phase 1.}
At step $t$, the process is in phase 1 if $X_t = Y_t = (l_t, m_t)$. Since both systems are in the same state and thus have the same service completion rate, there are two possible events for the coupled systems: 
\begin{enumerate}
    \item a rider joining the queue for both systems;
    \item a trip completion for both systems. 
\end{enumerate}
For both events, we have two possible subcases described as follows. First, if the number of dispatched drivers under $\phi^*$ and $\phi$ is the same, then $X_{t+1} = Y_{t+1}$, and the process remains in phase 1 at step $t+1$. 

Second, if the two policies dispatch a different number of drivers, the coupled states $X_{t+1} \neq Y_{t+1}$. Since the system will lie in the same anti-diagonal after dispatching, we know both systems lie in the same anti-diagonal at step $t+1$ under each event. Then, the process leaves phase 1 and enters phase 2 at step $t+1$.

\smallskip

\textbf{Phase 2.}
At step $t$, the process is in phase 2 if $X_t=(l_t,m_t)\neq Y_t=(l'_t,m'_t)$ and $l_t+m_t=l_t^\prime+m_t^\prime$ (i.e., $X_t$ and $Y_t$ are on the same anti-diagonal line). By Lemma \ref{lemma:serv_rate_ineq}, the trip completion rate follows $l_t\mu_{l_t,m_t}\geq l_t^\prime\mu_{l_t^\prime,m_t^\prime}$, which means system~1 has a larger trip completion rate than system~2 at the coupled state. Thus, it is possible that only system~1 experiences a trip completion at the next event. %
There are three possible events for the coupled systems: 
\begin{enumerate}
    \item a rider joining the queue for both systems;
    \item a trip completion for both systems;
    \item a trip completion only for system 1 (due to larger trip completion rate in system 1).
\end{enumerate}
For events (1) and (2), as the number of dispatched drivers does not affect the sum of $l_{t+1}$ and $m_{t+1}$, both systems will lie on the same anti-diagonal line at step $t+1$. If $X_{t+1}\neq Y_{t+1}$, the process stays at phase 2. If $X_{t+1}=Y_{t+1}$, the process leaves phase 2 and re-enters phase 1. \

For event (3), we have $l_{t+1}+m_{t+1}<l_{t+1}^\prime+m_{t+1}^\prime$, and the process leaves phase 2 and enters into phase 3.

\smallskip

\textbf{Phase 3.}
At step $t$, the process is in phase 3 if $l_t+m_t<l_t^\prime+m_t^\prime$. In this phase, we cannot determine which state, $X_t$ or $Y_t$, has the higher trip completion rate. Thus, there are four possible events for the coupled systems: 
\begin{enumerate}
    \item a rider joining the queue for both systems;
    \item a trip completion for both systems;
    \item a trip completion only for system 1 (if system 1 has a larger trip completion rate); 
    \item a trip completion only for system 2 (if system 2 has a larger trip completion rate).
\end{enumerate}
When events (1), (2), and (3) occur, the process stays at phase 3 since one can easily show that $l_{t+1}+m_{t+1}<l_{t+1}^\prime+m_{t+1}^\prime$ holds.

When event (4) occurs, there are two possible sub-cases: $l_{t+1}^\prime+m_{t+1}^\prime= l_{t+1}+m_{t+1}$ and $l_{t+1}^\prime+m_{t+1}^\prime>l_{t+1}+m_{t+1}$. If $l_{t+1}^\prime+m_{t+1}^\prime= l_{t+1}+m_{t+1}$, then the process leaves phase 3 and re-enters either phase 1 (when $X_{t+1}=Y_{t+1}$) or phase 2 (when $X_{t+1}\neq Y_{t+1}$).
Otherwise, if $l_{t+1}^\prime+m_{t+1}^\prime>l_{t+1}+m_{t+1}$, the process stays at phase 3.

\smallskip

In all phases above, since the effective request rate for both systems are kept the same, the revenue rate for both systems are the same as well. However, under coupling, the penalty rate in system 1 is consistently less or equal to that in system 2 because of weakly smaller $l+m$ values. Hence, we have $\Tilde{\mathcal{R}}(\phi^*,\lambda^*)\geq\Tilde{\mathcal{R}}(\phi^*,\lambda^{\mathrm{D}})\ge \Tilde{\mathcal{R}}(\phi,\lambda)$, which proves the result. The first inequality holds because there always exists an optimal Markovian policy (see Theorem 8.1.2 of \citealt{Puterman1994MDP}).

\medskip

\noindent\underline{\textbf{Case 2: $c_d>c_r$}}. 
The proof idea here is similar to the first case. We show that there exists a policy $(\phi^D,\lambda^D)$ such that $\Tilde{\mathcal{R}}(\phi^{\mathrm{D}},\lambda^{\mathrm{D}})\geq \Tilde{\mathcal{R}} (\phi,\lambda)$ holds for any policy $(\phi,\lambda)$. Construct $\phi^{\mathrm{D}}$ such that  $\phi^{\mathrm{D}}(l,m)=0$ for all type 1 states $(l,m)$. %
We consider a coupling of two systems, where system 1 uses policy $(\phi^{\mathrm{D}},\lambda^{\mathrm{D}})$ and system 2 uses policy $(\phi,\lambda)$. Similar to the case of $c_d=c_r$, $\lambda^{\mathrm{D}}$ here is defined to keep the effective request rate in system 1 always the same as that in system 2. Again, 
let \(\{X_t\}_{t=1}^\infty\) and \(\{Y_t\}_{t=1}^\infty\) denote the coupled Markov chain, where each step $t$ corresponds to a transition in the joint process. 
In the following, we describe the coupling processes and the constructed dispatching policy $\phi^{\mathrm{D}}$ in each phase.

\smallskip

\textbf{Phase 1.} At step $t$, the process is in phase 1 if $X_t = Y_t$. Let $X_t=(l_t, m_t)$ denote the corresponding state of system 1 at step $t$. There are two possible events for the coupled systems, which are the same as those in case 1, phase 1. In this phase, we set $\phi^{\mathrm{D}}(l,m)=\phi(l,m)$, which is the dispatching policy in system 2, for all type 2 states $(l,m)$. We now consider next state $Y_{t+1}=(l_{t+1}^\prime,m_{t+1}^\prime)$ for system 2. Let $P^\ast$ be the corresponding zigzag path under $\phi^\ast$, the optimal dispatching policy when $c_d=c_r$. There are two sub-cases for $Y_{t+1}$:
\begin{itemize}
    \item $Y_{t+1}$ is a type 2 state or $Y_{t+1}\in P^\ast$. Then by the construction of $\phi^{\mathrm{D}}$, both systems dispatch the same number of drivers at step $t$, so we have $X_{t+1} = Y_{t+1}$. The process remains in phase 1. 
    \item $Y_{t+1}$ is a type 1 state and $Y_{t+1}\notin P^\ast$. the number of dispatched drivers differs between the two policies, so we have $X_{t+1} \neq Y_{t+1}$. Since $\phi^{\mathrm{D}}(l,m)=0$ for all type 1 states $(l,m)$, we know $X_{t+1}\in P^\ast$. Moreover, $X_{t+1}$ and $Y_{t+1}$ still lie on the same anti-diagonal, and $l_{t+1}<l^\prime_{t+1}$. The process leaves phase 1 and enters into phase 2.  
\end{itemize}

\smallskip

\textbf{Phase 2.}
At step $t$, the process is in phase 2 if $l_{t}<l^\prime_{t}$, $m_t>m_t^{\prime}$, $l_t+m_t=l_t^\prime+m_t^\prime$, and $X_t\in P^\ast$. In other words, $Y_t$ lies on the same anti-diagonal as $X_t$ but is positioned to the lower left of $X_t$. In this phase, we set $\phi^{\mathrm{D}}=\phi^\ast$, which is the optimal dispatching policy when $c_d=c_r$. The coupled systems face the same three events defined in phase 2, case 1 ($c_d=c_r$). Again, we use the same label (1), (2), and (3) described in phase 2, case 1 to represent each event. 

When event (1) or (2) occurs, if $X_{t+1}\neq Y_{t+1}$, we know $l_{t+1}<l^\prime_{t+1}$ as $\phi^\ast$ makes at most one dispatch after an event occurs. Thus, the process stays at phase 2. If $X_{t+1}=Y_{t+1}$, the process leaves phase 2 and re-enters into phase 1. 

When event (3) occurs, we have $l_{t+1}<l_{t+1}^\prime$ and $l_{t+1}+m_{t+1}<l'_{t+1}+m'_{t+1}$, and the process leaves phase 2 and enters phase 3.

\smallskip

\textbf{Phase 3.}
At step $t$, the process is in phase 3 if $l_t<l_t^\prime$ and $l_t+m_t<l_t^\prime+m_t^\prime$. There are four possible events for the coupled systems, which are exactly the same as those in phase 3, case 1. We use event (1), (2), (3), and (4) to describe each event. We set dispatching policy $\phi^{\mathrm{D}}$ to be dependent on the event that occurred. 

When event (1), (2), or (3) occur, the two systems keep their anti-diagonal relationship, that is, $l_{t+1}+m_{t+1}<l_{t+1}^\prime+m_{t+1}^\prime$. Under these events, we set $\phi^{\mathrm{D}}(l,m)= 0$ for all type 2 states $(l,m)$, which means system 1 will not make a dispatch under these events. Then $l_{t+1}<l_{t+1}^\prime$ and $l_{t+1}+m_{t+1}<l_{t+1}^\prime+m_{t+1}^\prime$ hold due to the following facts:
\begin{itemize}
    \item For event (1), we have $m_{t+1}=m_t+1>m_t^\prime+1\geq m_{t+1}^\prime$ and $l_{t+1}=l_{t}<l_t^\prime \leq l_{t+1}^\prime$.
    \item For events (2) and (3), we have $m_{t+1}=m_t>m_t^\prime \geq m_{t+1}^\prime$ and $l_{t+1}=l_{t}-1<l_t^\prime-1\leq l_{t+1}^\prime$. 
\end{itemize} 
Thus, the process will stay at phase 3 since all conditions for phase 3 are satisfied. 

\smallskip

When event (4) occurs, we have $l_{t+1}+m_{t+1}=l_t+m_t$ and $l_{t+1}^\prime+m_{t+1}^\prime=l_t^\prime+m_t^\prime-1$. We construct $\phi^{\mathrm{D}}$ according to the following two sub-cases.
\begin{itemize}
    \item $l_{t+1}+m_{t+1}<l_{t+1}^\prime+m_{t+1}^\prime$. We set $\phi^{\mathrm{D}}(l,m) = 0$ for all type 2 states $(l,m)$ as well. Then we know $m_{t+1}=m_t>m_t^\prime\geq m_{t+1}^\prime$, and thus $l_{t+1}<l_{t+1}^\prime$ holds. As a result, the process will stay at phase 3. 
    \item $l_t+m_t=l_{t+1}^\prime+m_{t+1}^\prime$. We set $\phi^{\mathrm{D}}=\phi^\ast$, the optimal dispatching policy when $c_d=c_r$. It is easy to verify that the process will enter either phase 1 if $l_{t+1}=l_{t+1}^\prime$ or phase 2 if $l_{t+1}<l_{t+1}^\prime$.
\end{itemize}

\smallskip
 
In all phases above, the revenue rates are identical for both systems, while the penalty rate in system 1 is always less than or equal to that in system 2 at coupled states because system 1 has weakly smaller $l+m$ and $l$ values. Hence, again the objective function satisfies $\Tilde{\mathcal{R}}(\phi^*,\lambda^*)\geq\Tilde{\mathcal{R}}(\phi^{\mathrm{D}},\lambda^{\mathrm{D}})\geq \Tilde{\mathcal{R}}(\phi,\lambda)$, which proves the result.

\medskip

\noindent\underline{\textbf{Case 3: $c_d<c_r$}}. Similar to case 2, we show there exists $\phi^{\mathrm{D}}$ and $\lambda^{\mathrm{D}}$ such that $\Tilde{\mathcal{R}}(\phi^{\mathrm{D}},\lambda^{\mathrm{D}})\geq \Tilde{\mathcal{R}} (\phi,\lambda)$ for an arbitrary policy $(\phi,\lambda)$. We set $\lambda^{\mathrm{D}}$ to keep the effective request rate in system 1 always the same as that in system 2. We set $\phi^{\mathrm{D}}$ such that  $\phi^{\mathrm{D}}(l,m)=1$ for all type 2 states $(l,m)$. For type 1 states, $\phi^{\mathrm{D}}$ depends on system 2, and will be defined later. %
Again, let \(\{X_t\}_{t=1}^\infty\) and \(\{Y_t\}_{t=1}^\infty\) denote the coupled Markov chain, where each step $t$ corresponds to a transition in the joint process. By the partial construction of $\phi^\text{D}$, we know $X_t$ is a type 1 state for all $t$. There are four possible phases in this case, which are illustrated below.

\smallskip

\textbf{Phase 1.} At step $t$, the process is in phase 1 if $X_t = Y_t$. In this phase, there are two possible cases for $Y_{t+1}$ and the construction of $\phi^{\mathrm{D}}$ depends on whether $Y_{t+1}$ is a type 1 state or type 2 state.
\begin{itemize}
    \item $Y_{t+1}$ is a type 1 state. We set $\phi^{\mathrm{D}}=\phi$ for all type 1 states. Then the next coupled states satisfy $X_{t+1} = Y_{t+1}$, and the process remains in phase 1.
    \item $Y_{t+1}$ is a type 2 state. We set $\phi^{\mathrm{D}}=0$ for all type 1 states (i.e., $\phi^{\mathrm{D}}=\phi^\ast$). The next coupled states satisfy $X_{t+1} \neq Y_{t+1}$ and $X_{t+1}\in P^\ast$ since $\phi^{\mathrm{D}}(l,m)=1$ for all type 2 states. Since the state always stay on the same anti-diagonal after a dispatching action, we have $l_{t+1}+m_{t+1}=l_{t+1}^\prime+m_{t+1}^\prime$ and $m_{t+1}<m^\prime_{t+1}$. Then the process leaves phase 1 and enters phase 2.
\end{itemize}

\smallskip

\textbf{Phase 2.} At step $t$, the process is in phase 2 if $m_t<m_t^\prime$, $l_t>l_t^\prime$, $l_t+m_t=l_t^\prime+m_t^\prime$, and $X_t\in P^\ast$. Thus, $Y_t$ must be a type 2 state. There are three possible events for the coupled systems, which are the same as those in phase 2, case 1. Again, there are two possible cases for $Y_{t+1}$ and the construction of $\phi^{\mathrm{D}}$ depends on whether $Y_{t+1}$ is a type 1 state or type 2 state.
\begin{itemize}
    \item $Y_{t+1}$ is a type 1 state.  We set $\phi^{\mathrm{D}}=\phi$. Only event (1) or (2) can occur here because $Y_{t}$ is a type 2 state and $Y_{t+1}$ is a type 1 state, so there must be a state transition in system 2. Thus, we have $l_{t+1}+m_{t+1}=l^\prime_{t+1}+m^\prime_{t+1}$. Moreover, we know $X_{t+1}=Y_{t+1}$ because both systems have the same dispatching policy at type 1 states and $Y_{t+1}$ becomes a type 1 state now. Consequently, the process leaves phase 2 and enters phase 1 again.
    \item $Y_{t+1}$ is a type 2 state. We set $\phi^{\mathrm{D}}=\phi^\ast$, the optimal dispatching policy under $c_d=c_r$. Then $X_{t+1}\in P^\ast$ is a type 1 state. When event (1) or (2) occurs, we have $l_{t+1}+m_{t+1}=l^\prime_{t+1}+m^\prime_{t+1}$ because dispatching does not change the anti-diagonal of the state. Since $Y_{t+1}$ is a type 2 state and $X_{t+1}$ must be a type 1 state, we know $m_{t+1}<m_{t+1}^\prime$ and $l_{t+1}>l_{t+1}^\prime$ hold. Thus, the process stays in phase 2. 
    
    When event (3) occurs, we have $l_{t+1}+m_{t+1}=l_t+m_t-1<l^\prime_{t}+m^\prime_{t}=l^\prime_{t+1}+m^\prime_{t+1}$
    Moreover, $m_{t+1}$ either equals $m_t$ (not dispatch) or $m_t-1$ (dispatch), which is strictly less than $m_{t}^\prime$. The process then enters phase 3. 
\end{itemize}

\smallskip

\textbf{Phase 3.} At step $t$, the process is in phase 3 if $m_t<m_t^\prime$, $l_t+m_t<l_t^\prime+m_t^\prime$, and $X_t\in P^\ast$. There are four possible events for the coupled systems, which are exactly the same as those in phase 3, case 1. 
We construct $\phi^{\mathrm{D}}$ based on whichever of $m_t$ and $m_{t+1}^\prime$ is larger.

If $m_{t}>m_{t+1}^\prime$, at all type 1 states, we set $\phi^{\mathrm{D}}$ to keep dispatching until $m_{t+1}=m_{t+1}^\prime$. Note that event (3) is not possible to occur under $m_t>m_{t+1}^\prime$. In specific, we set $\phi^{\mathrm{D}}$  according to the following under events (1), (2), and (4): 

\smallskip

\begin{enumerate}
    \item[(a)] Under event (1), for every $k=0,\cdots,\min\{L-l_t,m_t-m_{t+1}^\prime+1\}$, we set $\phi^{\mathrm{D}}(l_t+k,m_t+1-k)=1$;
    \item[(b)] Under event (2), for every $k=0,\cdots,\min\{L-l_t+1,m_t-m_{t+1}^\prime\}$, we set $\phi^{\mathrm{D}}(l_t-1+k,m_t-k)=1$;
    \item[(c)] Under event (4), for every $k=0,\cdots,\min\{L-l_t,m_t-m_{t+1}^\prime\}$, we set $\phi^{\mathrm{D}}(l_t-1+k,m_t-k)=1$.
\end{enumerate}
For any type 1 state $(l,m)$ not covered above, we set $\phi^{\mathrm{D}}(l,m)=0$. 

On the other hand, if $m_t\leq m_{t+1}^\prime$, we set $\phi^{\mathrm{D}}= \phi^*$, the optimal dispatching policy under $c_d = c_r$. 

\smallskip

Now we discuss what happens for systems 1 and 2 under each event. When events (1) and (2) occur, there are three possible cases, which are $m_t> m_{t+1}^\prime$, $m_t= m_{t+1}^\prime$, and $m_t< m_{t+1}^\prime$.

\begin{itemize}
    \item If $m_t> m_{t+1}^\prime$ due to a dispatching in system 2, we make dispatches in system 1 to state $m_{t+1}=m_{t+1}^\prime$. Since both systems have an arrival or a trip completion, we have $l_{t+1}+m_{t+1}<l_{t+1}^\prime+m_{t+1}^\prime$. Then, the process enters into phase 4.
    \item If $m_{t+1}=m_{t+1}^\prime$, we have $\phi^{\mathrm{D}}= \phi^*$, so we know $X_{t+1}\in P^\ast$ and $l_{t+1}+m_{t+1}<l_{t+1}^\prime+m_{t+1}^\prime$. Thus, the process enters phase 4.
    \item If $m_{t+1}=m_{t+1}^\prime$, we have $\phi^{\mathrm{D}}= \phi^*$, the process stays in phase 3. 
\end{itemize}

When event (3) occurs, the state $Y_{t+1}$ is the same as $Y_t$, so we have $m_t<m_t^\prime=m_{t+1}^\prime$. As a result, $\phi^{\mathrm{D}}$ is set as $\phi^*$. Consequently, $m_{t+1}\leq m_t<m_t^\prime=m_{t+1}^\prime$ and  $X_{t+1}\in P^\ast$ hold. The process stays at phase 3. When event (4) occurs, we have $l_t+m_t=l_{t+1}+m_{t+1}$. There are two sub-cases, depending on whether the two systems are on the same anti-diagonal at $t+1$.
\begin{itemize}
    \item $l_t+m_t=l_{t+1}+m_{t+1}=l_{t+1}^\prime+m_{t+1}^\prime$, that is, the two systems are on the same anti-diagonal. Since $m_{t+1}\leq m_t$ holds, we know $m_{t+1}\leq  m_{t+1}^\prime$ holds according to our constructed $\phi^{\mathrm{D}}$. If $m_{t+1}^\prime=m_{t+1}$, the process re-enters phase 1. If $m_{t+1}<m_{t+1}^\prime$, we know $\phi^{\mathrm{D}}$ is constructed as $\phi^\ast$. Thus, the process re-enters phase 2.
    \item $l_t+m_t=l_{t+1}+m_{t+1}<l_{t+1}^\prime+m_{t+1}^\prime$. If $m_t<m_{t+1}^\prime$, we set $\phi^{\mathrm{D}}=\phi^\ast$, so the state in system 1 does not change (i.e., $l_{t+1}=l_t$ and $m_{t+1}=m_t$). Thus, the process stays at phase 3. If $m_t> m_{t+1}^\prime$, according to constructed $\phi^{\mathrm{D}}$, we have $m_{t+1}=m_{t+1}^\prime$. The process enters phase 4. If $m_t= m_{t+1}^\prime$, we set $\phi^{\mathrm{D}}=\phi^\ast$, and the process also enters phase 4 due to $m_{t+1}=m_t=m_{t+1}^\prime$.
\end{itemize}

\smallskip

 \textbf{Phase 4.} At step $t$, the process is in phase 4 if $m_t=m_t^\prime$, $l_t<l_t^\prime$, and $l_t+m_t<l_t^\prime+m_t^\prime$ hold. There are four possible events for the coupled system, which are exactly the same as those in phase 3, case 1. We construct $\phi^{\mathrm{D}}$ according to the state type of $Y_{t+1}$. If $Y_{t+1}$ is a type 1 state, we use the same construction of $\phi^{\mathrm{D}}$ as one in phase 3 under $m_t>m_{t+1}^\prime$. That is, at all type 1 states, $\phi^{\mathrm{D}}$ is set to keep dispatching until $m_{t+1}\leq m_{t+1}^\prime$ (i.e., dispatching the same number of drivers in both systems if system 1 is still in a type 1 state before taking the action).
Otherwise, if $Y_{t+1}$ is a type 2 state, we set $\phi^{\mathrm{D}}=\phi^*$. 

Now we discuss how states change under each event. When events (1) and (2) occur, if $Y_{t+1}$ is a type 1 state, we
have two subcases: 
\begin{itemize}
    \item $(l_{t},m_{t}+1)$ is a type 1 state under event (1) or $(l_{t}-1,m_{t})$ is a type 1 state under event (2). In this case, both systems dispatch the same number of drivers. Thus, we have $m_{t+1}=m_{t+1}^\prime$ and the process stays in phase 4.
    \item $(l_{t},m_{t}+1)$ is a type 2 state under event (1) or $(l_{t}-1,m_{t})$ is a type 2 state under event (2). If system 2 dispatches at least one driver, then according to the construction of $\phi^{\mathrm{D}}$, both systems dispatch the same number of drivers as well, so the process stays in phase 4. If system 2 does not dispatch any drivers, system 1 will dispatch exactly one driver since we require $\phi^{\mathrm{D}}(l,m)=1$ for all type 2 states. Therefore, it follows $m_{t+1}=m_t-1=m_t^\prime-1=m_{t+1}^\prime-1<m_{t+1}^\prime$ and $X_{t+1}\in P^\ast$ after dispatching, so the process enters phase 3. 
\end{itemize}

Otherwise, if $Y_{t+1}$ is a type 2 state under event (1) or (2), there is no dispatching in system 2 at step $t$ since $Y_t$ is a type 1 state. By the zigzag property, we know both $(l_{t}-1,m_{t})$ and $(l_{t},m_{t}+1)$ are type 2 states, so system 1 will make a dispatch. Then we have $m_{t+1}<m_{t+1}^\prime$, and $X_{t+1}\in P^\ast$. Thus, the process re-enters phase 3. 

When event (3) occurs, $Y_{t+1}$ is the same as $Y_t$, so we have $m_{t+1}^\prime=m_{t}^\prime=m_t$. If $(l_t-1,m_t)$ is a type 1 state, we have $m_{t+1}=m_t=m_t^\prime=m_{t+1}^\prime$, so the process stays at phase 4. On the other hand, if $(l_t-1,m_t)$ is a type 2 state, system 1 will make a dispatch and go into a type 1 state. Thus, we have $m_{t+1}<m_t=m_t^\prime=m_{t+1}^\prime$ and $X_{t+1}\in P^\ast$, so the process re-enters phase 3. 

When event (4) occurs, $Y_{t+1}$ must be a type 1 state. This is because $X_t=(l_t,m_t)$ is a type 1 state, and $l_t\leq l_{t+1}^\prime$ and $m_t\geq m_{t+1}^\prime$ hold, so by Lemma \ref{lemma:serv_rate_ineq}, $Y_{t+1}=(l_{t+1}^\prime,m_{t+1}^\prime)$ is also a type 1 state. Then, according to the construction of $\phi^{\mathrm{D}}$, we have $m_{t+1}=m_{t+1}^\prime$ since the number of dispatches under the event is the same in both system 1 and system 2. Thus, if $l_t+m_t<l_{t+1}^\prime+m_{t+1}^\prime$ holds, the process stays at phase 4. Otherwise, if $l_t+m_t=l_{t+1}^\prime+m_{t+1}^\prime$, the next states $X_{t+1}=Y_{t+1}$, so the process returns to phase 1.

In all phases, the revenue rate is identical for both systems, while the penalty rate in system 1 is always less than or equal to that in system 2 at coupled states because system 1 has weakly smaller $l+m$ and $m$ values. Hence, again the objective function satisfies $\Tilde{\mathcal{R}}(\phi^*,\lambda^*)\geq\Tilde{\mathcal{R}}(\phi^{\mathrm{D}},\lambda^{\mathrm{D}})\geq \Tilde{\mathcal{R}}(\phi,\lambda)$, which proves the result.\hfill\halmos

\end{proof}

\medskip

\begin{proof}{Proof of Theorem 2.}
By Theorem \ref{theorem:opt_policy_struct}, when $c_d=c_r$, there exists an optimal dispatching policy $\phi^\ast$ that satisfies $\phi^\ast(l,m)=0$ for all type 1 states and $\phi^\ast(l,m)=1$ for all type 2 states. Let $P^\ast=\big((l_1,m_1),(l_2,m_2),\cdots\big)$ be the corresponding zigzag path under $\phi^\ast$ and  $P_{j}^\ast:=\big((l_1,m_1),\cdots,(l_j,m_j)\big)$ be a sub-path of $P^\ast$ that terminates at state $(l_j,m_j)\in P^\ast$. In the following, we show that when $c_d=c_r$, $P_{j}^\ast$ yields an objective based on optimal static pricing that is greater than or equal to any other zigzag paths terminating at state $(l_j,m_j)$ when we compare two paths in Algorithm \ref{alg:DP_zigzag} (i.e., $P_{l_j,m_j}=P_j^\ast$ for $j=1,2,\cdots$, where $P_{l,m}$ stores the best zigzag path terminated at $(l, m)$ when Algorithm \ref{alg:DP_zigzag} terminates). 

We use induction to show this. Recall that in Algorithm~\ref{alg:DP_zigzag}, we evaluate the objective values of the two predecessor paths leading into state $(l_j, m_j)$: one arriving from above and the other from the left. For $k=1,2,\cdots$, suppose that path $P_{k}^\ast$ follows $P_{l_k,m_k}=P_k^\ast$. This holds trivially for $k=1$. We want to show $P_{l_{k+1},m_{k+1}}=P_{k+1}^\ast$ holds as well. We have the following two cases.

\smallskip

\begin{enumerate}
    \item[(1)] $l_{k+1}=l_{k}$ and $m_{k+1}=m_k+1$, that is, $P_{\text{left}}=P_{k+1}^\ast$ when Algorithm \ref{alg:DP_zigzag} is at $(l_{k+1},m_{k+1})$. We show that $\mathcal{R}_{\text{left}}\geq \mathcal{R}_{\text{above}}$. Similar to the proof idea of Theorem \ref{theorem:opt_policy_struct}, we construct a coupling between two systems: system 1 follows path $P_{\text{left}}$ and system 2 follows $P_{\text{above}}$. Both systems adopt the same static pricing, and we align their cutoff states so that the cutoff state in system 1 lies on the same anti-diagonal (same $l+m$ value) as the cutoff state in system 2. Under this coupling, the revenue rates are identical when neither system is at its cutoff state. Moreover, for states at the same anti-diagonal, the trip completion rate in system 1 is greater than or equal to that in system 2, so its penalty rate is always no greater than that of system 2 under coupling. In addition, system 1 never reaches its cutoff state if system 2 is not at its cutoff state. Together, these facts imply that the objective of system 1 is at least as large as that of system 2, i.e., $\mathcal{R}_{\text{left}}\geq\mathcal{R}_{\text{above}}$. If the objectives under the two paths satisfy $\mathcal{R}_{\text{left}}>\mathcal{R}_{\text{above}}$, according to Algorithm \ref{alg:DP_zigzag}, we have $P_{l_{k+1},m_{k+1}}=P_{k+1}^\ast$. Otherwise, if we have $R_{\text{left}}=R_{\text{above}}$, since $(l_{k+1}-1,m_{k+1})$ must be a type 2 state as it is not on the optimal zigzag path $P^\ast$, according to Algorithm \ref{alg:DP_zigzag}, we have $P_{l_{k+1},m_{k+1}}=P_{k+1}^\ast$ as well.
    
    \item[(2)] $l_{k+1}=l_{k}+1$ and $m_{k+1}=m_k$, that is, $P_{\text{above}}=P_{k+1}^\ast$. Under a similar coupling argument, one can show that $\mathcal{R}_{\text{above}}\geq\mathcal{R}_{\text{left}}$. If the objectives satisfy $\mathcal{R}_{\text{above}}>\mathcal{R}_{\text{left}}$, we have $P_{l_{k+1},m_{k+1}}=P_{k+1}^\ast$. If the objectives satisfy $\mathcal{R}_{\text{left}}=\mathcal{R}_{\text{above}}$, then $(l_{k+1}-1,m_{k+1})$ must be a type 1 state. Thus, according to Algorithm \ref{alg:DP_zigzag}, we have $P_{k+1}^\ast=P_{l_{k+1},m_{k+1}}$.
\end{enumerate}

\smallskip

Therefore, by induction, Algorithm \ref{alg:DP_zigzag} will set $P_{l_j,m_j}=P^\ast$ for every $j=1,2,\cdots$. Thus, %
Algorithm 1 yields the same optimal dispatching policy described in Theorem \ref{theorem:opt_policy_struct}. Moreover, the corresponding pricing policy $\lambda^\ast$ converges to the optimal solution under value iteration, so the pricing policy is optimal as well. \hfill\halmos  

\end{proof}

\medskip

\begin{proof}{Proof of Proposition \ref{prop:M_limit}.}
    Consider the threshold $M \geq L \bar{\mu} (p_0+p_{\mathrm{max}}t_0)/c_r$. First, we notice that for any $l \in \mathcal{L}$ and any $m \geq M$, we have
    \be \label{eq:p_max_less_penalty}
    p_0+p_{\mathrm{max}}t_0 \leq \frac{c_d l + c_r m}{l \mu_{l,m}}.
    \ee
    To see why, we start with the chosen threshold condition:
    \be
    M \geq \frac{L \bar{\mu} (p_0+p_{\mathrm{max}}t_0)}{c_r} \implies \frac{c_r M}{L \bar{\mu}} \geq (p_0+p_{\mathrm{max}}t_0).
    \ee
    Since $m \geq M$ and $\bar{\mu} \geq \mu_{l,m} > 0$, we have
    \be
    \frac{c_dl+c_r m}{l \mu_{l,m}} \geq \frac{c_r M}{l \bar{\mu}} \geq \frac{c_r M}{L \bar{\mu}} \geq (p_0+p_{\mathrm{max}}t_0)
    \ee
    for any $l \leq L$. Now, consider extending the path $P=\big((l_1,m_1),\cdots,(l_I,m_I)\big)$ to some state $(l_{I+1}, m_{I+1})$ where $m_{I+1} \geq M$. Let the extended path be $P^\prime:=\big((l_1,m_1),\cdots,(l_{I+1},m_{I+1})\big)$. We show that the objective under the optimal dynamic pricing for $P^\prime$ is the same as that for path $P$. To illustrate this, we examine a cycle of recurrence of the process that starts and ends at $(l_I, m_I)$. Within this cycle, if the first event is a trip completion, then the expected return in this cycle will be the same as that of the original path $P$ starting at $(l_I, m_I)$. On the other hand, if the first event is an effective arrival, the state transitions to $(l_{I+1}, m_{I+1})$. Then the expected time for the Markov chain to return to state $(l_I, m_I)$ is $1/(l_{I+1} \mu_{l_{I+1},m_{I+1}})$. During this period, the expected accumulated penalty is    $(c_dl_{I+1}+c_rm_{I+1})/(l_{I+1} \mu_{l_{I+1},m_{I+1}})$. However, the expected revenue coming from a single arrival is at most $p_0+p_{\mathrm{max}}t_0$, which is less than or equal to the penalty by equation \eqref{eq:p_max_less_penalty}. Thus optimal dynamic pricing on $P'$ will not admit any new rider arrival at state $(l_I, m_I)$.

    Clearly, the same argument works for optimal static pricing under path comparison in Algorithm \ref{alg:DP_zigzag}: when $m > M$, extending the path from state $(l-1,m)$ or $(l,m-1)$ to state $(l,m)$ will not improve the objective under any static pricing. Therefore, for $m>M$ and $l\in\mathcal{L}$, the objective recorded at each state $R_{l,m}$ satisfies $R_{l,m}=\max\{R_{l,m-1},R_{l-1,m}\}$. As a consequence, for $m>M$, the following equation holds:
    \begin{align*}
        R_{L,m}
        &=\max\{R_{L,m-1},R_{L-1,m}\}\\
        &=\max\{R_{L,m-1},R_{L-1,m-1},R_{L-2,m}\}\\
        &=\max\{R_{L,m-1},R_{L-1,m-1},\cdots,R_{1,m}\}\\
        &=\max_{l\in\{1,\cdots,L\}}\{R_{l,m-1}\}\\
        &=\max_{l\in\{1,\cdots,L\}}\{R_{l,M}\}=R_{L,M}.
    \end{align*}

    Thus, extending a path to any state $(l,m), m>M$ cannot produce a better solution. %
    \hfill\halmos 
\end{proof}

\medskip

\begin{proof}{Proof of Proposition \ref{prop:simulation_mu_property}.}
Given the estimated parameters $\hat{\alpha}_1, \hat{\alpha}_2$ and $\hat{C}$, 
the predicted service rate at state $(l,m)$ can be written as
\begin{align*}
    \hat{\mu}_{l,m}=\frac{1}{\hat{C}(m+1)^{\hat{\alpha}_2}(L-l+1)^{\hat{\alpha}_1}+t_0}.
\end{align*}

Since the fitted coefficient $\hat{\alpha}_1,\hat{\alpha}_2<0$, we know $\hat{\mu}_{l,m}$ is strictly increasing in $m$ and strictly decreasing in $l$. To verify Assumption \ref{as:2}, we need to show that the differences in the trip completion rate, $l\mu_{l,m+1} - l\mu_{l,m}$ and $l\mu_{l,m} - (l+1)\mu_{l+1,m}$, are both non-increasing in $m$ and non-decreasing in $l$. Instead of directly proving this, we establish a stronger result: the partial differences $\mu_{l,m+1} - \mu_{l,m}$ and $\mu_{l,m} - \mu_{l+1,m}$ are non-increasing in $m$ and non-decreasing in $l$. We first perform a variable transformation on $l$ by replacing it with $o = L - l$, where $o$ represents the number of idle drivers. Under this transformation, with a slight abuse of notation, we denote by the service rate $\hat{\mu}(o,m)$. We further extend the domain of $o$ and $m$ to continuous values, with $o \in [0, L]$ and $m \in [0, +\infty]$. In specific, $\hat{\mu}(o, m)$ can be written as
\begin{align*}
    \hat{\mu}(o,m)=\frac{1}{\hat{C}(m+1)^{\hat{\alpha}_2}(o+1)^{\hat{\alpha}_1}+t_0}.
\end{align*}
Since $\hat{\mu}(o, m)$ is continuously differentiable, it suffices to show that its Hessian matrix is element-wise non-positive. More specifically, we need to confirm that:
\begin{align*}
    \frac{\partial^2 \hat{\mu}}{\partial o \partial m} \leq 0, \quad \frac{\partial^2 \hat{\mu}}{\partial o^2} \leq 0, \quad \text{and} \quad \frac{\partial^2 \hat{\mu}}{\partial m^2} \leq 0.
\end{align*}
Due to the symmetry in $o$ and $m$, it is enough to verify $\frac{\partial^2 \hat{\mu}}{\partial o^2} \leq 0$ and $\frac{\partial^2 \hat{\mu}}{\partial o \partial m} \leq 0$.
\be \label{eq:partial_om}
\frac{\partial^2 \hat{\mu}}{\partial o \partial m}=-\frac{\hat{C} \hat{\alpha}_2 \hat{\alpha}_1\left(t_0-\hat{C}(o+1)^{\hat{\alpha}_1}(m+1)^{\hat{\alpha}_2}\right)(m+1)^{\hat{\alpha}_2-1}(o+1)^{\hat{\alpha}_1-1}}{\left(t_0+\hat{C}(m+1)^{\hat{\alpha}_2}(o+1)^{\hat{\alpha}_1}\right)^3},
\ee
\be \label{eq:partial_oo}
\frac{\partial^2 \hat{\mu}}{\partial o^2}=-\frac{\hat{C} \hat{\alpha}_1\left(t_0(\hat{\alpha}_1-1)-\hat{C}(1+\hat{\alpha}_1)(m+1)^{\hat{\alpha}_2}(o+1)^{\hat{\alpha}_1}\right)(m+1)^{\hat{\alpha}_2}(o+1)^{{\hat{\alpha}_1-2}}}{\left(t_0+\hat{C}(m+1)^{\hat{\alpha}_2}(o+1)^{\hat{\alpha}_1}\right)^3}.
\ee
Since $\hat{\alpha}_1 < 0$, $\hat{\alpha}_2 < 0$, and $t_0 - \hat{C}(o+1)^{\hat{\alpha}_1}(m+1)^{\hat{\alpha}_2} \geq t_0 - \hat{C} \geq 0$, it follows from equation \eqref{eq:partial_om} that
$
\frac{\partial^2 \hat{\mu}}{\partial o\partial m} \leq 0.
$
Next, for the second derivative with respect to $o$, we use the fact that $\hat{\alpha}_1 < 0$ and we have:
\begin{align*}
& t_0(\hat{\alpha}_1-1)-\hat{C}(1+\hat{\alpha}_1)(m+1)^{\hat{\alpha}_2}(o+1)^{\hat{\alpha}_1}\\
=& -t_0-\hat{C}(m+1)^{\hat{\alpha}_2}(o+1)^{\hat{\alpha}_1}+\hat{\alpha}_1\left(t_0-\hat{C}(m+1)^{\hat{\alpha}_2}(o+1)^{\hat{\alpha}_1}\right)\\
\leq & \hat{\alpha}_1\left(t_0-\hat{C}(m+1)^{\hat{\alpha}_2}(o+1)^{\hat{\alpha}_1}\right) \\
\leq & 0,
\end{align*}
which implies that $\frac{\partial^2 \hat{\mu}}{\partial o^2}\leq 0$ from equation \eqref{eq:partial_oo}. This completes the proof.\hfill\halmos  %
\end{proof}

\end{document}